\documentclass[lettersize,journal]{IEEEtran}
\usepackage{amsmath,amsfonts,amsthm}
\usepackage{algorithmic}
\usepackage{algorithm}
\usepackage{array}
\usepackage[caption=false,font=normalsize,labelfont=sf,textfont=sf]{subfig}
\usepackage{textcomp}
\usepackage{stfloats}
\usepackage{url}
\usepackage{verbatim}
\usepackage{graphicx}
\usepackage{cite}
\usepackage{multirow}

\allowdisplaybreaks

\usepackage[dvipsnames]{xcolor}

\newtheorem{theorem}{Theorem}

\newtheorem{definition}{Definition}
\newtheorem{corollary}{Corollary}

\newcounter{boxlblcounter}  

\newenvironment{boxlabel}
  {\begin{list}
    {\arabic{boxlblcounter}}
    {\usecounter{boxlblcounter}
     \setlength{\labelsep}{0.25em}
     \setlength{\itemindent}{0em} 
     \setlength{\leftmargin}{0.15cm}
     \setlength{\rightmargin}{0 cm}
     
    }
  }
{\end{list}}

\begin{document}

\title{Intercept Function and Quantity Bidding in Two-stage Electricity Market with Market Power Mitigation}

\author{Rajni Kant Bansal, Yue Chen, Pengcheng You, Enrique Mallada
\thanks{R. K. Bansal and E. Mallada are with the Johns Hopkins University, Baltimore, MD 21218, USA {(email: \{rbansal3,mallada\}@jhu.edu)}}

\thanks{Y. Chen is with The Chinese University of Hong Kong, Shatin, NT, Hong Kong SAR, China {(email:yuechen@mae.cuhk.edu.hk)}}

\thanks{P. You is with the College of Engineering, Peking University, Beijing, China  {(email:pcyou@pku.edu.cn)}}

}



\maketitle

\begin{abstract}
Electricity markets typically operate in two stages, day-ahead and real-time. Despite best efforts striving efficiency, evidence of price manipulation has called for system-level market power mitigation (MPM) initiatives that substitute noncompetitive bids with default bids. Implementing these policies with a limited understanding of participant behavior may lead to unintended economic losses. In this paper, we model the competition between generators and inelastic loads in a two-stage market with stage-wise MPM policies. The loss of Nash equilibrium and lack of guarantee of stable market outcome in the case of conventional supply function bidding motivates the use of an alternative market mechanism where generators bid an intercept function. A Nash equilibrium analysis for a day-ahead MPM policy leads to a Stackelberg-Nash game with loads exercising market power at the expense of generators. A comparison of the resulting equilibrium with the standard market (not implementing any MPM policy) shows that a day-ahead policy completely mitigates the market power of generators. On the other hand, the real-time MPM policy increases demand allocation to real-time, contrary to current market practice with most electricity trades in the day-ahead market. Numerical studies illustrate the impact of the slope of the intercept function on the standard market.

\end{abstract}


\section{Introduction}



\IEEEPARstart{I}{n} the US, electricity markets typically achieve a supply-demand balance in two stages,  day-ahead and real-time~\cite{two_stage_design}. Although designed to allocate resources efficiently and prevent speculation, in practice, there are efficiency losses
indicated through price differences between the two stages~\cite{low_lmp_spread,market_power_strat_gen_adam}. The discrepancies arise,
for example, due to various combinations of uncertainty in load forecast, unscheduled
maintenance, and other factors that may include price manipulation by strategic
participants who seek to benefit from the market. To encourage market competition, several independent system operators (ISOs) use different mitigation strategies, often triggered locally at predefined market conditions like congestion, to identify and mitigate non-competitive offers in either of the stages~\cite{local_mpm_stanford}. 
Despite the execution of such local policies, some ISOs have documented intervals with non-competitive participant behavior (e.g., $\sim 2\%$ hours in the California ISO region~\cite{caiso_scop}), which has subsequently led to system-level market power mitigation (MPM) policy initiatives~\cite{caiso_proposal,caiso_extda-bundle1}. Such system-level policies planned separately for each stage seek to identify non-competitive bids and substitute them with default bids to address the market power concerns. 

These approximate default bids estimate generator operation costs based on the operator's knowledge of technology, fuel prices, operational constraints~\cite[\S 30.11.2.3]{caiso_tarif}, etc. 
However, the competition between participants in a market with a system-level policy may result in unexpected market outcomes. Our previous work~\cite{bansal2023market} on market equilibrium with supply function bidding highlights the lack of guarantee of stable market outcome. That is, when participants are strategic, it is not possible to guarantee the existence of a Nash equilibrium. This motivates the search for alternative mechanisms that can provide guarantees of the existence of an equilibrium and, thus, better mitigate market power. 

In this paper, we consider the use of intercept function bidding~{\cite{yue_prosumer}} as an alternative market participation strategy that provides several benefits from the standpoint of market power mitigation. Precisely, we model the competition between generators and loads in a two-stage settlement electricity market where each generator bids the intercept of a supply function seeking to maximize their aggregate profit. Meanwhile, loads bid demand quantities and seek to minimize their payment in the market. We study the competition among these participants and consider the effect of super-imposing default bids on the equilibrium outcome.
Since the market operator can estimate the generation cost~\cite{caiso_proposal}, we ideally model the execution of the default bid MPM policy at each stage of the market,  by substituting generator bids within such stage with the true their cost.~\cite{bansal_e_energy}. A real-time MPM policy results in a Nash game between generator and load participants. However, a day-ahead MPM policy leads to a Stackelberg-Nash game with loads as leaders and generators as followers. Our equilibrium analysis of the proposed market mechanism highlights several advantages over previous works. 


\vspace{.5ex}
\noindent
\textbf{Contributions:} The main contributions in this paper are summarized below.
\begin{boxlabel}
    \item \emph{Standard Market: }We model and study the competition between generators (bidding the intercept of a supply function) and loads (bidding quantity) on a market without any market power mitigation strategy (called here standard market). We show that the competitive equilibrium (assuming price-taking participants) is efficient, in the sense that it aligns with a hypothetical social planner. 
    Notably, as seen in \cite{pyou_discovery} for supply function bidding, demand does not favor any particular stage for placing its bids. 
    A Nash equilibrium analysis for strategic participants shows generators profiting from price manipulations at the expense of higher load payments.
    
    \item \emph{Real-time MPM:} We then consider applying the default-bid MPM policy in the real-time stage. We show that our intercept-bidding mechanism always ensures the existence of both a competitive equilibrium (in the price-taking case) and a Nash equilibrium (in the strategic case). The competitive equilibrium shares the same characteristics as the standard market. However, the Nash equilibrium results in the undesirable market outcome where all the demand clears in the real-time market, a situation highly undesired by system operators.

    \item \emph{Day-ahead MPM:} Finally, we analyze the impact of a day-ahead default-bid MPM policy. We show that though the competitive equilibrium once again shares the same properties as the standard competitive equilibrium, the competition of strategic participants leads to a Stackelberg-Nash {\cite{stack_nash_tnn}} game between loads and generators wherein loads collectively act as leaders. Our analysis shows that not only a Stackelberg-Nash equilibrium always exists, but it also has the intended effect of mitigating the generator's market power by shifting part of the benefits of price manipulation back to the loads.


\end{boxlabel}

\vspace{.5ex}
\noindent
\textbf{Related Work:} There exists vast literature modeling the competition between participants in electricity markets and identifying price manipulation opportunities in different market settings, e.g., in single-stage settlement markets~\cite{DR_strat_nekouei,DA_strat_atzeni} or a two-stage settlement markets~\cite {pyou_discovery,nathan_two_stage}, in energy markets~\cite{Consumer_energy_trading,anderson_mpm} or capacity markets~\cite{RTN_capacity_market,Cap_expansion_game}, and with perfect competition of (price-taking) participants~\cite{bansal2021market,Joh2004Efficiency} or imperfect competition of strategic participants~\cite{Li_Na_linear_SFE,bansal2023market}. 
Beyond understanding market power, the design of mitigation strategies has been an important subject of study for operators and academics. This includes, e.g., CAISO's local market power mitigation policies~\cite{local_mpm_stanford}, congestion penalties~\cite{congestion_tariff_mpm}, virtual transactions~\cite{virtual_transactions_mpm}, \textcolor{black}{ forward contracting~\cite{allaz_vila,cai_forward_contract}, demand shifting~\cite{yujian_ds_mpm}, capacity regulation~\cite{newbery2002mitigating}, etc.} 

However, the impact of system-level mitigation policies, such as those proposed by CAISO,  has received limited attention. To the best of our knowledge, our previous work~\cite{bansal2023market,bansal_e_energy} is the first such attempt to do a counter-factual analysis of the effect of default-bid system-level MPM policy on the market outcome and its participants. Particularly,~\cite{bansal2023market} considers such a setting for generators bidding conventional supply functions but cannot guarantee the existence of the Nash equilibrium. In~\cite{bansal_e_energy}, the use of intercept function bids is considered to model a day-ahead default-bid MPM policy. However, the analysis is limited as it does not compare the market outcome (Nash Equilibrium) that one would have obtained without such a policy. 
Building on \cite{bansal_e_energy}, in this paper, we provide a thorough counter-factual analysis of the effect of default-bid MPM policies both in day ahead and real time, providing convincing proof that intercept-bidding when combined with a default-bid day-ahead MPM policy constitutes an effective strategy that outperforms, both standard intercept bid market, and real-time default-bid MPM policies.

\vspace{.5ex}
\noindent
\textbf{Paper Organization:} 
The rest of the paper is structured as follows. In Section~\ref{sec_2}, we formulate the social planner problem, describe the two-stage market mechanism, and define participants’ behavior and two-stage market equilibrium. In Section~\ref{sec_3}, we characterize the market equilibrium in a standard two-stage market model based on intercept bidding. We model system-level market power mitigation policies and characterize the market equilibrium for different participation behavior in Section~\ref{sec_4}. We provide insights on the market outcome in a market with MPM policy and compare it with the standard market in Section~\ref{sec_5}. To streamline the presentation, we relegate the comparison of the intercept bid-based standard market with the slope bid-based standard market to Section~\ref{slope_compare_sec}. Finally, conclusions are in Section~\ref{sec_6}.

\subsubsection*{Notation} 
\textcolor{black}{
The standard notation $f(x,y)$ denotes a function of independent variables $x$ and $y$. We use $f(x;y)$ to represent a function of independent variable $x$ and parameter $y$. Also, $|\mathcal{I}|$ represents the cardinality of the set $\mathcal{I}$. 
}

\section{Market Model}\label{sec_2}

In this section, we start with the formulation of the underlying social planner problem. We then describe the standard two-stage settlement electricity market design, where a generator bids an intercept function while demand bids quantities, and define rational participants’ behavior, i.e., price-taking or price-anticipating. Finally, we define a general market equilibrium in such a market setting.

\subsection{Social Planner Problem}

Consider a single-interval two-stage settlement electricity market where a set $\mathcal{G}$ of generators compete with set $\mathcal{L}$ of inelastic loads. The power dispatch of generator $j$ over the two stages is denoted by $g_j \in  \mathbb{R}$ such that 
\begin{align}\label{gen_two_stage}
    g_j := g_j^{d} + g_j^{r}
\end{align}
where $g_j^d\in  \mathbb{R}, \ g_j^r\in  \mathbb{R}$ denote the dispatch in the two stages of day-ahead and real-time markets, respectively.  The total inelastic demand of load $l$ allocated over two stages is denoted by $d_l\in  \mathbb{R}^{+}$, such that 
\begin{align} \label{load_two_stage}
    d_l := d_l^d + d_l^r
\end{align}
where $d_l^d\in  \mathbb{R}, \ d_l^r\in  \mathbb{R}$ denote the allocated load in the day-ahead and real-time markets, respectively. The market operator seeks to achieve supply-demand balance, i.e. 
\begin{align} \label{power_bal}
    \sum\nolimits_{j\in\mathcal{G}}g_j = \sum\nolimits_{l \in \mathcal{L}}d_l
\end{align}
The social planner problem that seeks to minimize the cost of dispatching generators to meet aggregate demand is given by:
%
\begin{align}\label{planner_problem}
    \min_{g_j,j\in\mathcal{G}} & \ \sum\nolimits_{j\in \mathcal{G}} \frac{c_j}{2} g_j^2  \ \  \text{s.t.} \ \eqref{power_bal}
\end{align}
where we assume a quadratic cost of dispatching generators, parameterized by quadratic coefficients $c_j \in  \mathbb{R}^+$. 
The underlying social planner problem~\eqref{planner_problem} is considered a benchmark, and we will analyze the deviation between market equilibrium and the social planner solution as one of the metrics to study the existence of market power. 

\subsection{Two-stage Market Mechanism}
We now describe the two-stage market clearing.
\subsubsection{Day-ahead Market}
 Each generator $j$ submits an intercept function, with constant slope $b^d \in  \mathbb{R}_{>0}$ and parameterized by $\beta_j^d \in  \mathbb{R}$, that indicates the willingness of the generator to participate in the market, given by:
\begin{align}\label{gen_intercept_bid_da}
    g_j^d = b^d\lambda^d - \beta_j^d
\end{align}
where $\lambda^d$ denotes the day-ahead price. Each load $l \in  \mathcal{L}$ in the day-ahead market bids quantity $d_l^d$. Once all the bids $(\beta_j^d,d_l^d)$ are received, the market clears with supply-demand balance:   
\begin{align}\label{da_power_bal}
    \sum\nolimits_{j \in \mathcal{G}} \left(b^d\lambda^d - \beta_j^d\right) = \sum\nolimits_{l \in \mathcal{L}} d_l^d.
\end{align}

The solution to~\eqref{da_power_bal} gives the dispatch and clearing price such that generator $j$ earns $\lambda^d g_j^d$ while load $l$ pays $\lambda^d d_l^d$ in the market settlement process.

\subsubsection{Real-time Market}
 Similar to the day-ahead market, each generator $j$ submits an intercept function, with constant slope $b^r\in  \mathbb{R}_{>0}$ and parameterized by $\beta_j^r\in  \mathbb{R}$, as:
\begin{align}\label{gen_intercept_bid_rt}
    g_j^r = b^r\lambda^r - \beta_j^r
\end{align}
where $\lambda^r$ denotes the real-time prices. Each load $l \in  \mathcal{L}$ in real-time market bids quantity $d_l^r$. The load allocation in the real-time market is given once the load allocation in the day-ahead market is determined due to the demand inelasticity and \eqref{load_two_stage}. Once all the bids $(\beta_j^r,d_l^r)$ are received, the market clears with supply-demand balance, given by   
\begin{align}\label{rt_power_bal}
    \sum\nolimits_{j \in \mathcal{G}} \left(b^r\lambda^r - \beta_j^r\right) = \sum\nolimits_{l \in \mathcal{L}} d_l^r.
\end{align}

The solution to~\eqref{rt_power_bal} determines the dispatch and clearing price such that generator $j$ earns $\lambda^r g_j^r$ while load $l$ pays $\lambda^r d_l^r$ in the market settlement process.

\subsection{Participant Behavior}

We focus on two different forms of rational participation behavior, i.e., price-taking and price-anticipating, where each generator $j$ (load $l$) seeks to maximize (minimize) its profit (payment) in the two-stage market. The profit of generator $j$, denoted by $\pi_j$, is given by:
\begin{align}\label{generator_profit}
    &\pi_j(g_j^{d},g_j^{r},\lambda^{d},\lambda^{r}) \!:= \! \lambda^{r}g_j^{r} \!\!+ \lambda^{d}g_j^{d} \!-\! \frac{c_j}{2} (g_j^d+g_j^r)^2
\end{align}
Similarly, the payment of load $l$, denoted by $\rho_l$, is given by:
\begin{align}\label{load_payment_definition}
\!\!\rho_l(d_l^{d},d_l^{r},\lambda^{d},\lambda^{r}) & :=  \lambda^{d}d_l^{d} \!+\! \lambda^{r}d_l^{r} = \lambda^{d}d_l^{d} \!+\! \lambda^{r}(d_l \!-\! d_l^d) 
\end{align}
where we substitute the load inelaticity constraint~\eqref{load_two_stage}.


\subsubsection{Price-taking Participation} We first define the price-taking participant behavior and then formulate the individual problem of participants. 
\begin{definition}    
    A participant is price-taking in the market if it does not anticipate the impact of its bid on the market prices and accepts the existing prices as given. 
\end{definition}

Given the day-ahead and real-time prices $(\lambda^d,\lambda^r)$ in the market, the individual problem of price-taking generator $j$ is:

\begin{align}\label{generator_price_taking_profit}
    &\max_{g_j^{d},g_j^{r}} \  \pi_j(g_j^{d},g_j^{r};\lambda^{d},\lambda^{r}) 
\end{align}
and the individual problem of price-taking load $l$ is given by:
\begin{align}\label{load_price_taking_payment}
    \min_{d_l^{d}} \  &  \rho_l(d_l^{d};\lambda^{d},\lambda^{r})
\end{align}

\subsubsection{Price-anticipating (Strategic) Participation}

We now define the price-anticipating participant behavior.
\begin{definition}
    A participant is price-anticipating (strategic) in the two-stage market if it can manipulate the prices by anticipating the impact of its bid and knowledge of other participants’ bids in two stages.
\end{definition}

Given load bids $d_l^{d}, d_l^{r}, l\in\mathcal{L}$, and other generators' bids $\beta_k^d,\beta_k^r, k \in \mathcal{G}, k \neq j$, the individual problem of a price-anticipating generator $j$ is then given by:
\begin{subequations}\label{generator_strategic_profit_total}
\begin{align}\label{generator_strategic_profit}
    & \!\!\!\max_{g_j^{d},g_j^{r}} \ \! \pi_j \! \left(g_j^{d},g_j^{r},\lambda^{d}\!\left(g_j^{d};\overline{g}_{-j}^{d},d^{d}\right),\lambda^{r}\!\!\left(g_j^{r};\overline{g}_{-j}^{r},d^{r}\right) \right) \\
    &\textrm{ s.t. }  \ \ \ \ \eqref{da_power_bal}, \eqref{rt_power_bal}
\end{align}
\end{subequations}
where $\overline{g}_{-j}^{d} := \sum_{k \in \mathcal{G}, k \neq j}g_k^{d}$, and $\overline{g}_{-j}^{r} := \sum_{k \in \mathcal{G}, k \neq j}g_k^{r}$. Similarly, the individual problem for price-anticipating load $l$ is given by:
\begin{align}\label{load_strategic_payment}
    & \!\!\!\!\min_{d_l^{d}} \ \!\rho_l\!\left(\!d_l^{d}, \lambda^{d}\!\left(\!d_l^d;g_j^{d},\overline{d}_{-l}^{d}\!\right)\!\!,\lambda^{r}\!\!\left(\!d_l^d;g_j^{r},\overline{d}_{-l}^{r}\right)\!\right) \ \!\! \textrm{ s.t. } \eqref{da_power_bal}, \eqref{rt_power_bal}
\end{align}
where $\overline{d}_{-l}^{d} := \sum_{l \in \mathcal{L}, k \neq l}d_l^{d}, \ \overline{d}_{-l}^{r} := \sum_{l \in \mathcal{L}, k \neq l}d_l^{r}$.

\subsection{Market Equilibrium}

In this section, we describe the notion of market equilibrium 
in a two-stage settlement electricity market. In the market, firms make decisions in their best interest without accounting for others' incentives. However, at the equilibrium, the resulting prices are such that the market achieves the supply-demand balance, and no participating firm has any incentive to deviate from its bid. More formally,

\begin{definition} \label{market_eqbm}
A two-stage market is at equilibrium if the participant bids and market clearing prices $( {\beta_j^{d}}, {\beta_j^{r}},j\in \mathcal{G}, d_l^{d}, d_l^{r}, l \in\mathcal{L},\lambda^{d}, \lambda^{r})$ in the day-ahead and real-time markets satisfy the following conditions:
\begin{enumerate}
    \item The bid $\beta_j^{d}, \beta_j^{r}$ of generator $j$ maximizes its profit.
    \item The allocation $d_l^{d}, d_l^{r}$ of load $l$  minimizes its payment.
    \item  The market clears with prices $\lambda^{d}$ given by~\eqref{da_power_bal} and $\lambda^{r}$ given by~\eqref{rt_power_bal}.
\end{enumerate}
\end{definition}
An equilibrium analysis of the market is often used to understand the efficiency and stability of a market mechanism. Though equilibrium is hard to attain in reality due to the dynamic nature of the market, descriptive and predictive equilibrium outcomes (if possible) provide intuition about the behavior of individual participants~\cite{eqbm_theory_2011}. We use equilibrium analysis in this paper to analyze the impact of system-level MPM policies on market outcomes. 
\begin{definition}\label{symm_market_eqbm}
A market equilibrium that satisfies the Definition~\ref{market_eqbm} is said to be symmetric on the generator side if all the generators are homogeneous and make the same decisions in both stages, i.e., $\beta_j^d:=\beta^d,  \ \beta_j^r:= \beta^r,  \ \forall j \in \mathcal{G}$.
\end{definition}

Supply function equilibrium is hard to analyze in the closed form, and prior works have considered similar simplifying assumptions to gain insights. In this work, we make this assumption to characterize the Nash equilibria and investigate the impact of MPM policies.

\section{Equilibrium in Standard Market}\label{sec_3}

In this section, we model the competition between generators and loads in a standard two-stage market without any mitigation policy. The participants bid in both day-ahead and real-time markets. 
We analyze such a game backward, starting from the real-time market, for the equilibrium path. The resulting equilibrium is regarded as a benchmark to determine the impact of the stage-wise system-level MPM policies later. 

\subsection{Competitive Equilibrium}

We first consider the case of price-taking participants in the market. We substitute~\eqref{gen_intercept_bid_da} and \eqref{gen_intercept_bid_rt} into \eqref{generator_price_taking_profit} to get the individual problem of generator $j$, given the prices $(\lambda^d,\lambda^r)$, as:
\begin{align}\label{gen_std_intercept_profit_comp}
    \!\!\!\!\!\!\max_{\beta_j^d,\beta_j^r} \ \!\!\!-\beta_j^d{\lambda\!}^d \!-\!\beta_j^r{\lambda\!}^r \!\!-\! \frac{c_j}{2}(\!\beta_j^d\!+\!\beta_j^r\!)^2 \!\!+\! c_j(\!b^d{\lambda\!}^d\!+\!b^r{\lambda\!}^r\!)(\!\beta_j^d\!+\!\beta_j^r\!)    
\end{align}
The individual problem of load $l$ is given in the optimization problem~\eqref{load_price_taking_payment}. We can now characterize the competitive equilibrium in this market setting:

\begin{theorem}\label{comp_eqbm_wout_mpm_thrm} 
A competitive equilibrium in a standard two-stage settlement market without any mitigation policy exists and is given by
\begin{subequations}\label{comp_eqbm_wout_mpm} 
\begin{align}
    &{\beta}_j^d+{\beta}_j^r = \frac{b^d+b^r - c_j^{-1}}{\sum_{j\in\mathcal{G}}c_j^{-1}}d, \  \forall j \in \mathcal{G} \label{comp_eqbm_wout_mpm.a} \\
    & \sum_{j \in \mathcal{G}} (b^d \lambda^d \!-\! \beta_j^d) = \!\sum_{l \in \mathcal{L}}d_l^d, \ \sum_{j \in \mathcal{G}} (b^r \lambda^r \!-\! \beta_j^r) = \!\sum_{l \in \mathcal{L}}d_l^r \label{comp_eqbm_wout_mpm.b} \\
    & d_l^d+d_l^r = d_l, \forall l \in \mathcal{L} \label{comp_eqbm_wout_mpm.c} \\
    & \lambda^d = \lambda^r = \frac{1}{\sum_{j\in\mathcal{G}}c_j^{-1}}d \label{comp_eqbm_wout_mpm.d} 
\end{align}
\end{subequations}
\end{theorem}
We provide the proof of the theorem in Appendix~\ref{appendix_comp_eqbm_wout_mpm_thrm}. Although the competitive equilibrium in Theorem~\ref{comp_eqbm_wout_mpm_thrm} exists non-uniquely, i.e., each load $l$ is indifferent to demand allocation due to equal prices in the two stages, the resulting dispatch and prices align with the underlying social planner optimum~\eqref{planner_problem}. 

\subsection{Nash Equilibrium}

We next characterize the Nash equilibrium as a result of competition between price-anticipating participants. 
We first characterize the interaction between generators and loads in a real-time market for some given allocation in the day-ahead market. This results in a real-time subgame equilibrium that will help compute the Nash equilibrium in the two-stage market. 

\begin{theorem}\label{augmented_dispatch_thrm}
     We assume that there is more than one strategic generator in the market, i.e., $|\mathcal{G}| > 1$. The subgame equilibrium $(g_j^{r},d^r, \lambda^{r})$ due to the interplay between generators and loads in the real-time market, given the day-ahead market outcome $(g_j^{d},d_l^{d})$, is an optimal primal-dual solution to an augmented convex social planner problem, as:
    \begin{subequations}\label{augmented_dispatch}
    \begin{align}
        \min_{g_j^{r}} & \ \sum_{j\in\mathcal{G}}\left( \frac{1}{2b^r(|\mathcal{G}|-1)}{g_j^{r}}^2+\frac{c_j}{2}\left(g_j^{d}+g_j^{r}\right)^2\right) \label{augmented_obj}\\
        & \textrm{s.t.} \ \ \sum\nolimits_{j \in \mathcal{G}} g_j^r = \sum\nolimits_{l \in \mathcal{L}} d_l^r \label{augmented_constraint}
    \end{align}
    \end{subequations}
\end{theorem}
 We provide the proof of the theorem in Appendix~\ref{strat_eqbm_wout_mpm_augmented_proof}. The strategic participation of generators in real-time shifts the dispatch of generators, captured by the first term in the objective function of the augmented social planner problem in Theorem~\ref{augmented_dispatch_thrm}. Since the augmented problem is strictly convex, the subgame equilibrium is unique. Moreover, the subgame equilibrium does not exist if there is only one generator in the market and prices become indefinite. 

The following theorem characterizes the resulting symmetric Nash equilibrium in the market, where each individual generator solves \eqref{generator_strategic_profit_total} while each individual load solves \eqref{load_strategic_payment}. For tractability and closed-form analysis, we consider the participation of homogeneous generators in the market. 

\begin{theorem}\label{nash_eqbm_wout_mpm_thrm}
    Let's assume that generators are homogeneous, i.e., $c_j: = c,  \forall j \in \mathcal{G}$. If there is more than one generator participating in the market, i.e., $|\mathcal{G}|>1$, then the two-stage symmetric Nash equilibrium uniquely exists and it is given by:
    \begin{subequations}\label{nash_eqbm_wout_mpm}
    \begin{align}
          & \beta_j^d \!= \frac{b^dc}{|\mathcal{G}|}d + \frac{b^rc - \frac{|\mathcal{G}|-2}{|\mathcal{G}|-1}}{b^rc+\frac{|\mathcal{L}|+1}{|\mathcal{G}|-1}}\frac{|\mathcal{L}|+1}{|\mathcal{G}|(|\mathcal{G}|-1)}d^d, \ \forall j \in \mathcal{G} \label{nash_eqbm_wout_mpm.a1}\\
          & \beta_j^r \!=\!  \frac{b^rc}{|\mathcal{G}|}d \!-\! \frac{|\mathcal{G}|-2}{|\mathcal{G}|(|\mathcal{G}|-1)}d^r, \ \forall j \in \mathcal{G} \label{nash_eqbm_wout_mpm.a2}\\
          & g_j^d = \frac{1}{|\mathcal{G}|}d^d,  \ g_j^r = \frac{1}{|\mathcal{G}|}d^r, \ \forall j \in \mathcal{G}\label{nash_eqbm_wout_mpm.b}\\
          & \!\!d_l^d \!\!= \!\frac{b^dd_l}{b^d\!+\!b^r\!(\!|\mathcal{G}|\!\!-\!1\!)} \!+\! \frac{\frac{b^d}{1+b^rc(|\mathcal{G}|\!-\!1)}}{b^d\!+\!b^r\!(\!|\mathcal{G}|\!\!-\!1\!)}d^r \!\!\!-\! \frac{b^r}{b^d\!+\!b^r\!(\!|\mathcal{G}|\!\!-\!1\!)}d^d\label{nash_eqbm_wout_mpm.c1}\\
          & \! d_l^r = d_l - d_l^d, \forall l \in \mathcal{L}  \label{nash_eqbm_wout_mpm.c2}\\
          & \lambda^d \!= \!\frac{b^rc(|\mathcal{G}|\!-\!1) \!+\!2}{b^rc(|\mathcal{G}|\!-\!1) \!+\! 1}\frac{c}{|\mathcal{G}|}d+\frac{\left(\frac{b^r}{b^d}\!-\!1\right)c \!+\! \frac{1}{b^d(|\mathcal{G}|-1)}}{b^rc(|\mathcal{G}|\!-\!1) \!+\! 1}\frac{d^d}{|\mathcal{G}|}, \label{nash_eqbm_wout_mpm.d}\\
          & \!\!\lambda^r \!\!=\! \lambda^d \!+\!\frac{\frac{1}{|\mathcal{G}|(|\mathcal{G}|\!-\!1)}\left(\frac{|\mathcal{G}|\!-\!2}{|\mathcal{G}|\!-\!1} \!-\! b^rc\right)d}{b^d\!\left(\!b^rc\!+\!\frac{|\mathcal{L}|+1}{|\mathcal{G}|-1}\!\right) \!\!+\! b^r\!\!\left(\!b^rc\!+\!\frac{1}{|\mathcal{G}|-1}\!\right)\!(|\mathcal{G}|\!+\!|\mathcal{L}|\!\!-\!\!1\!)\!} \label{nash_eqbm_wout_mpm.e}
    \end{align}
    \end{subequations}
\end{theorem}
We provide the proof of the theorem in Appendix~\ref{app_nash_eqbm_wout_mpm}. At the equilibrium, the load allocation across stages depends on the slope of the bidding function, and operators can tune these for a higher allocation in the day-ahead market, as observed in current market practice. More specifically, we provide such a condition on the slope of the intercept functions in Corollary~\ref{load_corollary}. Moreover, for $|\mathcal{G}| = 1$, the generator makes arbitrary large bid decisions to drive prices high in the market, and the Nash equilibrium does not exist.



\begin{corollary}\label{load_corollary}
    \textcolor{black}{The load allocation across the two stages at the Nash equilibrium in a standard market~\eqref{nash_eqbm_wout_mpm} is given by}:
    \begin{subequations}
    \begin{eqnarray}
        & d^{d} = \frac{b^d\left(b^rc+\frac{|\mathcal{L}|+1}{|\mathcal{G}|-1}\right)}{b^d\left(b^rc+\frac{|\mathcal{L}|+1}{|\mathcal{G}|-1}\right) + b^r\left(b^rc+\frac{1}{|\mathcal{G}|-1}\right)(|\mathcal{G}|+|\mathcal{L}|-1)}d 
        \\ 
\label{load_allocation_day_ahead}
    & d^{r} = \frac{b^r\left(b^rc+\frac{1}{|\mathcal{G}|-1}\right)(|\mathcal{G}|+|\mathcal{L}|-1)}{b^d\left(b^rc+\frac{|\mathcal{L}|+1}{|\mathcal{G}|-1}\right) + b^r\left(b^rc+\frac{1}{|\mathcal{G}|-1}\right)(|\mathcal{G}|+|\mathcal{L}|-1)}d 
\label{load_allocation_real_time}
    \end{eqnarray}
    \end{subequations}
    Furthermore, for 
    \[
        b^d \ge b^r\frac{\left(b^rc+\frac{1}{|\mathcal{G}|-1}\right)(|\mathcal{G}|+|\mathcal{L}|-1)}{\left(b^rc+\frac{|\mathcal{L}|+1}{|\mathcal{G}|-1}\right)},
    \]
    the load allocation in the day-ahead market is higher than in the real-time market, i.e., $d^d \ge d^r$.
    
\end{corollary}

\section{Equilibrium in Market with an MPM Policy}\label{sec_4}

In this section, we model the impact of system-level MPM policies on market equilibrium. Each generator operates truthfully in the stage with an MPM policy in response to operator intervention in the form of a mitigation policy. With considerable market knowledge of participants' technology, fuel prices, operational constraints, historical prices, etc., ISOs can estimate, if not accurately, a reasonable bound on the operation cost of generators, which is used in substituting their bids with default bids in the presence of an MPM policy. However, each generator is allowed to bid an intercept function in the other stage. 

These policies are planned firstly for the real-time followed by the day-ahead market to keep a check on the high risk of market power exercise in the real-time market compared to the day-ahead market. For this paper, we assume that the operator can accurately estimate the operation cost of a generator such that the market clears efficiently in the stage with MPM policy. Such assumptions will allow us to develop an understanding of the system-level MPM policies and compare them with the standard market.

\subsection{Real-time MPM Policy}

In this subsection, we model the real-time default-bid MPM policy, formulate the individual problem for different participation behavior and characterize the market equilibrium.

\subsubsection{Modelling Real-time Default-bid MPM Policy}

For the real-time MPM policy, the operator accurately estimates the truthful operation of the generators in the real-time market, given the dispatch in the day-ahead market, i.e., 
\begin{align} \label{rt_true_dispatch}
        g_j^{r} = c_j^{-1}\lambda^{r}- g_j^{d} , \ \forall j \in  \mathcal{G}
\end{align}
Summing the equation~\eqref{rt_true_dispatch} over $j \in \mathcal{G}$ and substituting the two-stage supply-demand balance~\eqref{power_bal}, we get
\begin{align}\label{rt_true_prc}
    {\lambda^{r}} = \frac{d}{\sum_{j \in \mathcal{G}}c_j^{-1}}
\end{align}

We characterize the equilibrium in a market with a real-time MPM in the following subsection.

\subsubsection{Competitive Equilibrium}

We first consider the case of price-taking participants in the market. We substitute~\eqref{gen_intercept_bid_da}, \eqref{rt_true_dispatch}, and \eqref{rt_true_prc} in \eqref{generator_price_taking_profit} to get the individual problem of price-taking generator $j$, given the clearing price $\lambda^d$, as:
\begin{align}\label{generator_profit_rt_mpm}
    &\max_{\beta_j^d} \ \tilde{\pi}_j(\beta_j^{d};\lambda^{d}) \!:= \! \max_{\beta_j^d} \left(\frac{d}{\sum_{j \in \mathcal{G}}c_j^{-1}} - \lambda^{d}\right)\beta_j^d
\end{align}
Similarly, substituting~\eqref{rt_true_prc} in \eqref{load_price_taking_payment} gives the individual problem of load $l$ as:
\begin{align}\label{load_payment_intermediate_rt_mpm}
    & \min_{d_l^d}  \ \tilde{\rho}_l(d_l^{d};\lambda^{d}) := \min_{d_l^d} \ \left(\lambda^{d} - \frac{d}{\sum_jc_j^{-1}}\right)d_l^{d}
\end{align}
where the price $\lambda^d$ is given in the market. The resulting competitive equilibrium is characterized below:
\begin{theorem} \label{comp_eqbm_rt_mpm_thrm}
The competitive equilibrium in a two-stage market with a real-time MPM policy exists and it is given by:
\begin{subequations}\label{comp_eqbm_rt_mpm}
\begin{align}
    & g_j^{d} + g_j^{r} \!= \frac{c_j^{-1}}{\sum_{j\in\mathcal{G}}c_j^{-1}}d,  \ \beta_j^d \in \mathbb{R} \ \forall j\in\mathcal{G}\label{comp_eqbm_rt_mpm.a}\\
    & d_l^{d} + d_l^{r} = d_l, \ \!\! \forall l\in\mathcal{L}\label{comp_eqbm_rt_mpm.b}\\
    & \lambda^{d} = \lambda^{r} = \frac{1}{\sum_{j\in\mathcal{G}}c_j^{-1}}d\label{comp_eqbm_rt_mpm.c}
\end{align}
\label{competitive_eqbm_rt_mpm}
\end{subequations}
\end{theorem}
We provide proof of the theorem in Appendix~\ref{app_comp_eqbm_rt_mpm}. The competition between generators and loads for a higher and lower price market, respectively, leads to divergent behaviors. A set of equilibria exists in the market for equal prices in two stages. However, at such an equilibrium, loads do not have any incentive to allocate demand in the day-ahead market. Interestingly, the resulting competitive equilibrium still aligns with the social planner problem~\eqref{planner_problem}.

\subsubsection{Nash Equilibrium}

In this section, we characterize the market equilibrium for the competition between price-anticipating participants. Substituting~\eqref{rt_true_dispatch} and \eqref{rt_true_prc} in \eqref{generator_strategic_profit_total}, we get the individual problem of the price-anticipating generator $j$ that seeks to maximize the profit as:
\begin{subequations}\label{generator_strategic_profit_total_rt_mpm}
\begin{align}\label{generator_strategic_profit_rt_mpm}
    & \!\!\max_{\beta_j^{d},\lambda^d} \ \! \pi_j \! \left(\beta_j^{d},\lambda^{d}\!\left(\beta_j^{d};\overline{\beta}_{-j}^{d},d^{d}\right)\right) \\
    &\textrm{ s.t. } \eqref{da_power_bal}
\end{align}
\end{subequations}
Similarly, we substitute~\eqref{rt_true_dispatch},\eqref{rt_true_prc} in \eqref{load_strategic_payment} to get the individual problem of the price-anticipating load as:
\begin{subequations}\label{load_strategic_payment_rt_mpm}
\begin{align}
    & \!\!\!\min_{d_l^{d},\lambda^d} \ \rho_l\!\left(d_l^{d}, \lambda^{d}\!\left(d_l^d;\beta_j^{d},\overline{d}_{-l}^{d}\right)\right)\\
    &\textrm{ s.t. } \eqref{da_power_bal}.
\end{align}
\end{subequations}
We analyze the sequential game backward, starting with the real-time market where generators operate truthfully, resulting in fixed clearing prices. Although loads could bid in the real-time market, the bids are fixed by their decisions in the day-ahead market and load inelasticity. Therefore, each participant competes in the day-ahead market for individual interests. The following theorem characterizes the Nash equilibrium.

\begin{theorem}\label{strat_eqbm_rt_mpm_thrm}
If there is more than one generator participating in the market, i.e., $|\mathcal{G}|>1$, the two-stage Nash equilibrium in a market with a real-time MPM policy uniquely exists, as:
\begin{subequations}
\begin{align} \label{strat_eqbm_traditional_rt_mpm.b0}
    & g_j^{d} = 0,  \!\ g_j^{r} \!=\! \frac{c_j^{-1}}{\sum\nolimits_{j\in\mathcal{G}}c_j^{-1}}d, \ \beta_j^{d} \!=\! \frac{b^d}{\sum\nolimits_{j \in  \mathcal{G}}c_j^{-1}}d, \ \!\forall j\in\mathcal{G} \\ \label{strat_eqbm_traditional_rt_mpm.b1}
    & d_l^{d} = 0, \ d_l^{r} = d_l, \ \forall l\in\mathcal{L} \\ \label{strat_eqbm_traditional_rt_mpm.b2}
    & \lambda^{d} = \lambda^{r} = \frac{1}{\sum\nolimits_{j\in\mathcal{G}}c_j^{-1}}d   
\end{align}\label{strat_eqbm_traditional_rt_mpm}
\end{subequations}    
\end{theorem}
We provide proof of the theorem in Appendix~\ref{app_strat_eqbm_rt_mpm}. For a non-zero demand allocation in the day-ahead market, generators have the incentive to change their bid while attempting to manipulate prices and extract higher profit. Loads attempt to decrease prices to seek minimum payment simultaneously. The mutual competition to outbid each other results in the same price across stages, and all the demand shifts to the real-time market. Although there is no price difference across stages, i.e., no arbitrage opportunity, and the market dispatch aligns with the social planner optimum, i.e., efficient market equilibrium, such an equilibrium may not be desirable from the operator's perspective. In practice, day-ahead accounts for a majority of energy trades.

\subsection{Day-ahead MPM Policy}

In this section, we consider the impact of a day-ahead MPM policy. 

\subsubsection{Modeling Day-ahead Default-bid MPM policy} In this case, operators make an accurate estimation of generator dispatch cost in the day-ahead market, i.e.,
\begin{align} \label{da_true_dispatch}
    g_j^{d} =c_j^{-1}\lambda^d
\end{align}
Summing the equation~\eqref{da_true_dispatch} over $j \in \mathcal{G}$ and using the power-balance in day-ahead market~\eqref{da_power_bal} implies that:
\begin{align}\label{da_true_prc}
    \lambda^{d} = \frac{d^d}{\sum_{j\in\mathcal{G}}c_j^{-1}}
\end{align}

Each generator has the flexibility to bid in the real-time market and we characterize the resulting market equilibrium in the following subsection.

\subsubsection{Competitive Equilibrium}

We first define the individual problem of participants and then characterize the resulting competitive equilibrium. The individual problem of price-taking generator $j$ is given by:
\begin{align}
    \!\!\! \max_{\beta_j^r} \ \!\! \tilde{\pi}_j(\beta_j^{r};\!{\lambda\!}^{r}\!) \! := \!\max_{\beta_j^r} \ \! \! \! -\beta_j^r{\lambda\!}^{r} \!\!-\!  \frac{c_j}{2}\!\!{\left(\!\!\frac{c_j^{-1}d^d}{\sum\nolimits_{j \in \mathcal{G}}c_j^{-1}}\! +\! {b}^r{\lambda\!}^r \!\!-\! \beta_j^r\!\!\right)\!\!\!}^2 \label{generator_price_taking_profit_bids_da_mpm}
\end{align}
where we substitute~\eqref{da_true_dispatch},\eqref{da_true_prc} in \eqref{generator_price_taking_profit}. Similarly, the individual problem of load $l$ is given by~\eqref{load_price_taking_payment}. The resulting competitive equilibrium is characterized in the theorem below.
\begin{theorem} \label{comp_eqbm_da_mpm}
The competitive equilibrium in the two-stage market with a day-ahead MPM policy exists:
\begin{subequations}\label{comp_eqbm_da_mpm_equations}
\begin{align}
    & g_j^{d} = \frac{c_j^{-1}}{\sum\nolimits_{j\in\mathcal{G}}c_j^{-1}}d, \ g_j^{r} \!=\! 0, \ \beta_j^d \!=\! \frac{b^d}{\sum\nolimits_{j\in\mathcal{G}}c_j^{-1}}d, \ \forall j\in\mathcal{G} \label{comp_eqbm_da_mpm.a}\\
    & d_l^{d} + d_l^{r} = d_l \ \forall l\in\mathcal{L}, \ d^{d} = d, \ d^{r} = 0 \label{comp_eqbm_da_mpm.b}\\
    & \lambda^{d} = \lambda^{r} = \frac{1}{\sum\nolimits_{j\in\mathcal{G}}c_j^{-1}}d \label{comp_eqbm_da_mpm.c}
\end{align}
\label{competitive_eqbm_traditional}
\end{subequations}
\end{theorem}
The proof of the theorem was first presented in our previous paper~\cite{bansal_e_energy}, we include it here in Appendix~\ref{app_comp_eqbm_da_mpm} for completeness. Unlike the case of the real-time MPM policy in Theorem~\ref{comp_eqbm_rt_mpm_thrm} with equal prices across stages, the equilibrium in Theorem~\ref{comp_eqbm_da_mpm} is unique and incentivizes load to allocate all the demand in the day ahead market. It is efficient, i.e., aligns with the social planner problem, and is desirable from the current market practice perspective. 

\subsubsection{Nash Equilibrium}

We next consider the competition between price-anticipating participants in a market with a day-ahead MPM policy. The sequential game where generators operate truthfully in the day-ahead market results in a Stackelberg-Nash game with loads making decisions in the day-ahead as leaders and generators participating as followers in the real-time market. In addition to a leader-follower game between loads and generators, participants compete with each other within their own group for individual interests in a Nash game. Similar formulations have been analyzed in the literature for various market settings~\cite{stack_nash_lodi,stack_nash_tnn}. In this paper, we use the terminology from the reference~\cite{stack_nash_tnn}.  



The following theorem characterizes the Nash equilibrium, where load minimizes its payment as a leader, anticipating the prices in two stages and with knowledge of others' bids.  

\begin{theorem} \label{strat_eqbm_da_mpm_thrm}
Assume there is more than one generator participating in the market with a day-ahead MPM policy, i.e., $|\mathcal{G}|> 1$. Then the Nash equilibrium exists uniquely as:
\begin{subequations}\label{strat_eqbm_traditional_equations}
\begin{align}
    \label{strat_eqbm_traditional.a}
    & \!\! g_j^{d} \!= \!\! \left(\!\!1\!-\!\frac{\sum\limits_{j\in\mathcal{G}}\!\!{\left(\!\!\frac{1}{b^r(\!|\mathcal{G}|-1\!)}\!+\!c_j\!\!\right)\!\!}^{-1}}{(|\mathcal{L}|+1)\sum\limits_{j\in\mathcal{G}}\!\!c_j^{-1}}\!\right)\!\!\frac{c_j^{-1}}{\sum\limits_{j\in\mathcal{G}}\!\!c_j^{-1}}d, \ \forall j\in\mathcal{G}\\
    \label{strat_eqbm_traditional.b1}
    & g_j^{r} =\frac{{\left(\!\!\frac{1}{b^r(\!|\mathcal{G}|-1\!)}\!+\!c_j\!\!\right)\!\!}^{-1}}{|\mathcal{L}|+1}\frac{d}{\sum\limits_{j\in\mathcal{G}}c_j^{-1}}, \ \forall j\in\mathcal{G} \\
    \label{strat_eqbm_traditional.b2}
    & \!\!\!\!\!\! \beta_j^d \!\!= \!\!\!\left(\!\!b^d\!\!-\!\frac{b^d\!\!\sum\limits_{j\in\mathcal{G}}\!\!{\left(\!\!\frac{1}{b^r(\!|\mathcal{G}|-1\!)}\!+\!c_j\!\!\right)\!\!}^{-1}}{(|\mathcal{L}|+1)\!\!\sum\limits_{j\in\mathcal{G}}\!\!c_j^{-1}}\!-\!\frac{{\left(\!\!\frac{1}{b^r(\!|\mathcal{G}|-1\!)}\!+\!c_j\!\!\right)\!\!}^{-1}}{|\mathcal{L}|+1}\!\!\right)\!\!\frac{d}{\sum\limits_{j\in\mathcal{G}}\!\!c_j^{-1}}\!\\
    \label{strat_eqbm_traditional.e}    
    & d_l^{d} \!=\! d_l\!+\!\!\left(\!\frac{1}{|\mathcal{L}|+1}d\!-\!d_l\!\!\right)\frac{\sum\limits_{j\in\mathcal{G}}{\left(\!\!\frac{1}{b^r(\!|\mathcal{G}|-1\!)}\!+\!c_j\!\!\right)\!\!}^{-1}}{\sum\limits_{j\in\mathcal{G}}c_j^{-1}}, \ \forall l\in\mathcal{L} \\ 
\label{strat_eqbm_traditional.c}
    & d_l^{r} = \left(d_l - \frac{1}{|\mathcal{L}|+1}d\right)\frac{\sum\limits_{j\in\mathcal{G}}{\left(\!\!\frac{1}{b^r(\!|\mathcal{G}|-1\!)}\!+\!c_j\!\!\right)\!\!}^{-1}}{\sum\limits_{j\in\mathcal{G}}c_j^{-1}}, \ \forall l\in\mathcal{L} \\
\label{strat_eqbm_traditional.d}
    & \lambda^{d} \!=\!  \frac{d}{\sum\limits_{j\in\mathcal{G}}\!\!c_j^{-1}} -\frac{1}{|\mathcal{L}|\!+\!1}\frac{\sum\limits_{j\in\mathcal{G}}{\left(\!\!\frac{1}{b^r(\!|\mathcal{G}|-1\!)}\!+\!c_j\!\!\right)\!\!}^{-1}}{\sum_{j\in\mathcal{G}}c_j^{-1}}\frac{d}{\sum\limits_{j\in\mathcal{G}}c_j^{-1}}, \\
\label{strat_eqbm_traditional.f}    
    & \lambda^{r} \!=\! \lambda^d +\!\frac{1}{|\mathcal{L}|+1}\frac{d}{\sum\nolimits_{j\in\mathcal{G}}c_j^{-1}}  
\end{align}
\end{subequations}
\end{theorem}
The proof of the theorem was first presented in our previous paper~\cite{bansal_e_energy}, and we include it here in Appendix~\ref{app_strat_eqbm_da_mpm} for completeness.
Unlike the standard two-stage Nash equilibrium in Theorem~\ref{nash_eqbm_wout_mpm_thrm}, in the presence of a day-ahead MPM policy, the resulting Nash equilibrium always leads to higher prices in the real-time market. As generators operate truthfully in the day-ahead market, loads exploit this opportunity to allocate higher demand in the day-ahead market to seek minimum payment. Generators, with the flexibility to bid in the real-time market, attempt to manipulate and drive prices in the real-time market. The design of the day-ahead MPM policy puts generators in an inherent disadvantage position as followers in the market.

\begin{corollary}
At the Nash equilibrium~\eqref{strat_eqbm_traditional_equations} in a market with a day-ahead MPM policy, the load allocation in the day-ahead and the real-time market is given by:
    \begin{subequations}
    \begin{eqnarray}
        & d^d \!=\!\! \sum\limits_{l\in\mathcal{L}} d_l^{d} \!=\!  \!\left(\!1\!-\!\frac{1}{|\mathcal{L}|+1}\frac{\sum\limits_{j\in\mathcal{G}}\!\!{\left(\!\frac{1}{b^r(\!|\mathcal{G}|-1\!)}\!+\!c_j\!\right)\!}^{-1}}{\sum\limits_{j\in\mathcal{G}}c_j^{-1}}\right)d \in \left(\frac{d}{2}, d\right) \\ 
\label{load_allocation_day_ahead_int}
    & d^r = \sum\limits_{l\in\mathcal{L}}  d_l^{r} =  \frac{1}{|\mathcal{L}|+1}\frac{\sum\limits_{j\in\mathcal{G}}\!\!{\left(\!\frac{1}{b^r(\!|\mathcal{G}|-1\!)}\!+\!c_j\!\right)\!}^{-1}}{\sum\limits_{j\in\mathcal{G}}c_j^{-1}}d \in \left(0, \frac{d}{2}\right)
\label{load_allocation_real_time_int}
    \end{eqnarray}
    \end{subequations}
\end{corollary}
The proof uses the relation $b^r>0$ and sums up the individual load allocation at the Nash equilibrium~\eqref{strat_eqbm_traditional_equations}.

\section{Market Analysis}\label{sec_5}

In this section, we analyze the impact of system-level mitigation policies by comparing the resulting market equilibria with standard market equilibrium. 

\subsection{Equilibrium Insights on Stage-wise MPM Policies}

\begin{table}[!t]
  \caption{Competitive Equilibrium (CE) and Nash Equilibrium (NE) with a stage-wise MPM policy}
  \label{tab:eqbm_table}
  \centering
  \begin{tabular}{c|c|c}
    \hline
    Instance \!\!\!& Real-Time MPM & Day-Ahead MPM\\
    \hline
    \multirow{3}{*}{CE} & Non-unique equilibrium \!&\!Unique equilibrium\\
    & Solves social planner & Solves social planner\\
    & Arbitrary load allocation \! & All load in day-ahead \\
    \hline
    \multirow{4}{*}{NE} & Unique and efficient & Unique and non-efficient\\
    & All load in real-time &  majority load in day-ahead\\
    & undesirable to operator & desired market power mitigation\\
    \hline
\end{tabular}
\end{table}

We first discuss the case of the real-time MPM policy followed by the day-ahead MPM policy, as summarized in Table~\ref{tab:eqbm_table}. The mitigation policies in real time result in equal prices across stages, which is the same as the system marginal cost, and the market outcome aligns with the social planner optimum~\eqref{planner_problem} at both competitive~\eqref{comp_eqbm_rt_mpm} and Nash equilibrium~\eqref{strat_eqbm_traditional_rt_mpm}. However, the competitive equilibrium outcome fails to incentivize loads to allocate demand in the day-ahead market~\eqref{comp_eqbm_rt_mpm.b}. On the other hand, Nash equilibrium incentivizes loads to allocate demand to the real-time market entirely~\eqref{strat_eqbm_traditional_rt_mpm.b1}, making it undesirable from the operators' perspectives.

The day-ahead MPM policy also results in a unique competitive equilibrium~\eqref{comp_eqbm_da_mpm_equations} that aligns with the social planner optimum~\eqref{planner_problem} while incentivizing loads to allocate demand to the day-ahead~\eqref{comp_eqbm_da_mpm.b}. At the Nash equilibrium, the mitigation policy leads to generators participating as followers and limiting their market power. Generators participate strategically in real-time, inflating the prices above the system marginal cost~\eqref{strat_eqbm_traditional.f}. However, loads acting as leaders anticipate the real-time sub-game equilibrium and allocate more demand in the day-ahead market~\eqref{strat_eqbm_traditional.e}. Although a higher demand allocation in the day-ahead market increases the day-ahead clearing prices~\eqref{strat_eqbm_traditional.d}, it is still below the clearing prices in the real-time market~\eqref{strat_eqbm_traditional.f}. The loads are favored in the competition with a total payment at Nash equilibrium below the competitive equilibrium levels, as shown in row 1 of Table~\ref{tab:cost_profit_table}. 

\begin{table*}
  \caption{Comparison of normalized Nash equilibrium (normalized with competitive equilibrium) between a standard market and a day-ahead market policy market }
  \label{tab:cost_profit_table}
  \resizebox{\linewidth}{!}{\begin{tabular}{c|c|c|c}
    \hline
    Case \!\!\!\!& Social Cost & Generators Aggregate Profit & Loads Aggregate Payment\\
    \hline
    \vspace{1.5 pt}
    \!\!DA-MPM \!\!\!\!\!\!&\!\!\!\!\!\! $1\!+\!\frac{1}{\sum\limits_{j\in\mathcal{G}}\!\!\!c_j^{-1}}\frac{\Delta}{(|\mathcal{L}|+1)^2}\!$\!\!\!\!&\!\!\!\!$1\!-\!\frac{\sum_{j\in\mathcal{G}}\!{\left(\!\!\frac{1}{b^r(\!|\mathcal{G}|-1\!)}\!+c_j\!\!\right)\!\!}^{-1}}{\sum_{j\in\mathcal{G}}\!c_j^{-1}}\frac{2|\mathcal{L}|}{(|\mathcal{L}|+1)^2}\!-\!\frac{1}{\sum_{j\in\mathcal{G}}\! c_j^{-1}}\frac{\Delta}{(|\mathcal{L}|+1)^2}$\!\!\!\!\!&\!\!\!\!$1\!-\!\frac{\sum_{j\in\mathcal{G}}\!{\left(\!\!\frac{1}{b^r(\!|\mathcal{G}|-1\!)}\!+c_j\!\!\right)\!\!}^{-1}}{\sum_{j\in\mathcal{G}}\!c_j^{-1}}\!\frac{|\mathcal{L}|}{(|\mathcal{L}|+1)^2}$\!\!\\
    \!\!Standard \!\!\!\!\!& 1 &\!\!\!$1 \!+\! \frac{2}{b^rc(|\mathcal{G}|-1)+1}\frac{d^dd^r}{d^2}\!+\!\frac{2}{b^dc(|\mathcal{G}|-1)}\frac{(d^d)^2}{d^2} \!+\! \frac{2}{b^rc(|\mathcal{G}|-1)}\frac{(d^r)^2}{d^2} \nonumber $\!\!\!&\!\!\!$1 \!+\! \frac{1}{b^rc(|\mathcal{G}|-1)+1}\frac{d^dd^r}{d^2}\!+\!\frac{1}{b^dc(|\mathcal{G}|-1)}\frac{(d^d)^2}{d^2} \!+\! \frac{1}{b^rc(|\mathcal{G}|-1)}\frac{(d^r)^2}{d^2} \nonumber $\!\!\!\!\!\\
    \hline
    \multicolumn{4}{l}{\text{where } $\Delta := \sum_{j\in\mathcal{G}}{\frac{c_j}{{\left(\!\!\frac{1}{b^r(\!|\mathcal{G}|-1\!)}\!+c_j\!\!\right)\!}^2}}-\frac{{\left(\sum_{j\in\mathcal{G}}{\left(\!\!\frac{1}{b^r(\!|\mathcal{G}|-1\!)}\!+ c_j\!\!\right)\!\!}^{-1}\!\!\right)\!}^2}{\sum_{j\in\mathcal{G}}c_j^{-1}}$} \\
\end{tabular}}
\end{table*}

\begin{corollary}
 In a market with a day-ahead MPM policy, the total generator profit at the Nash equilibrium \eqref{strat_eqbm_traditional_equations} is always below the competitive equilibrium levels \eqref{competitive_eqbm_traditional}. 
\end{corollary}

From the market perspective, the social cost is higher at the Nash equilibrium~\eqref{strat_eqbm_traditional_equations} than the competitive equilibrium~\eqref{competitive_eqbm_traditional}, as shown in column 1 of Table~\ref{tab:cost_profit_table}.

\begin{corollary}
Assuming generators are homogeneous, i.e., $c_j=c,~ \forall j \in \mathcal{G}$, the social cost at the Nash equilibrium~\eqref{strat_eqbm_traditional_equations} is the same as the competitive equilibrium~\eqref{competitive_eqbm_traditional}.
\end{corollary}

The corollary uses the fact that for homogeneous generators $\Delta = 0$, as shown in Table~\ref{tab:cost_profit_table}. The term $\Delta$ is a non-linear function of the cost coefficients of generators and provides a quantitative measure of the heterogeneity in the system.

\subsection{Comparison of DA-MPM Policy with a Standard Market}

We next compare only the equilibrium for a day-ahead MPM policy with equilibria in a standard market, as the real-time MPM policy market equilibrium results in undesirable market outcomes. Unlike a set of competitive equilibria in a standard market~\eqref{comp_eqbm_wout_mpm}, the competitive equilibrium in the market with a day-ahead MPM policy is unique and incentivizes loads to allocate demand in the day-ahead market~\eqref{comp_eqbm_da_mpm_equations}.  

Interestingly at the Nash equilibrium in a market with a day-ahead MPM policy, clearing prices in real-time is always higher than in the day-ahead market~\eqref{strat_eqbm_traditional.f} due to the leader-follower structure and strategic participation of generators in real-time only. However, in the standard market, generators exploit the inelasticity of demand to manipulate the prices at Nash equilibrium in two stages resulting in higher day-ahead clearing prices~\eqref{nash_eqbm_wout_mpm.e} under certain conditions, i.e., the number of generators participating in the market and slope of the intercept function. We study the role of price-anticipating participants in a standard market and market with a day-ahead mitigation policy from the market and individual perspectives, i.e.,  social cost, generators' profit, and loads' payment in Tables~\ref{tab:cost_profit_table}.

For the sake of comparison between two market settings, 
we evaluate the Nash equilibrium with the assumption that generates are homogeneous and participate symmetrically in the market. Since generators are homogeneous, the market clears with the minimum cost of dispatch that equals the social planner cost, as shown in column 1 of Table~\ref{tab:cost_profit_table}.
We next look at the individual perspective to evaluate the properties of the Nash equilibrium. In the standard market, generators win the competition at the Nash equilibrium since they always earn a higher profit than the one achieved in the competitive equilibrium level, as shown in row 2 of Table~\ref{tab:cost_profit_table}. However, in the case of the day-ahead MPM policy, loads win the competition with lower payment at the Nash equilibrium than the competitive equilibrium, as shown in row 1 of Table~\ref{tab:cost_profit_table}. Although the day-ahead MPM policy does have the intended mitigation effect on the market power of generators, it results in loads exercising market power at the expense of generators.

\begin{figure*}
  \centering
  \includegraphics[width=0.95\linewidth]{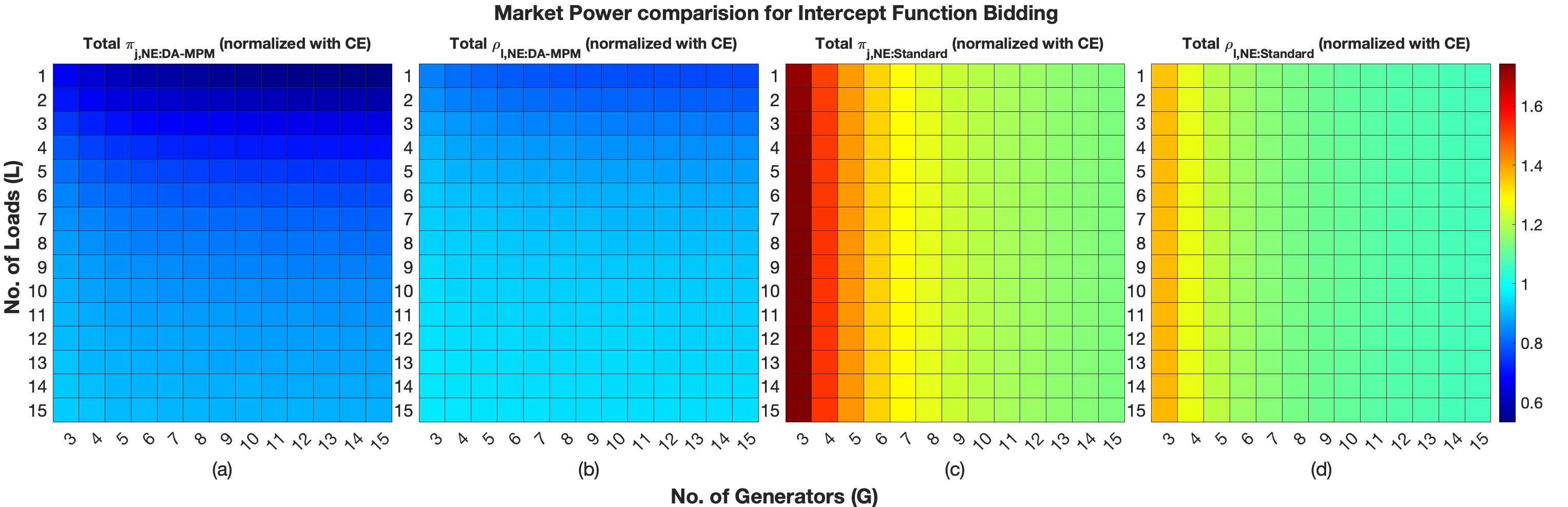}
  \caption{Total profit (a) and total payment (b) at Nash Equilibrium (NE) normalized with competitive equilibrium (CE) in a market with day-ahead MPM (DA-MPM), and total profit (c) and total payment (d) at Nash Equilibrium (NE) normalized with competitive equilibrium (CE) in a standard market.}
  \label{fig:ratio}
\end{figure*}

Figure~\ref{fig:ratio} compares the (normalized) aggregate profit and (normalized) aggregate payment at the Nash equilibrium in the standard market with a day-ahead MPM policy (DA-MPM) market, respectively. For simplicity, we assume that $b^d = b^r = \frac{1}{c}$. The aggregate generator profit (load payment) at the Nash equilibrium is normalized with the corresponding competitive equilibrium levels, which are the same in both market settings, and analyzed as we increase the number of participants in the market. The aggregate profit ratio in the DA-MPM policy market, as given by
\[
    1 - \frac{b^rc(\!|\mathcal{G}|-1\!)}{1+ b^rc(\!|\mathcal{G}|-1\!)}\frac{2|\mathcal{L}|}{(|\mathcal{L}|+1)^2} ,
\]
increases monotonically in the number of loads due to increased competition between loads, signaling a recovery in efficiency, whereas the ratio decreases monotonically in the number of generators due to increased competition between generators. This increased competition with an increase in the number of generators exacerbates their exploitation in the market, as shown by darker colors in the rows of panels (a) and (b) in Figure~\ref{fig:ratio}. 

However, the aggregate profit or payment ratio in the standard market increases with the number of loads and decreases with the number of generators, as shown in panels (c) and (d) in Figure~\ref{fig:ratio}.
This implies that generators always win the competition in the standard market with higher profit levels at the Nash equilibrium compared with the competitive equilibrium. Moreover, the day-ahead MPM policy results in the complete mitigation of generator market power, as shown in the comparison of generator normalized aggregate profit in the two markets in panels (a) and (c) in Figure~\ref{fig:ratio}, respectively.

\section{Equilibrium Comparison with Slope Function Bid in a Standard Market}\label{slope_compare_sec}

In this section, we compare the intercept function bidding with the conventional slope function bidding\footnote{For ease of comparison between the two bidding mechanisms, we say a generator submits an intercept function~\cite{yue_prosumer} or a slope function~\cite{pyou_discovery} when it bids intercept or slope of the supply function, respectively}, a.k.a. linear supply function in a standard market (without the implementation of an MPM policy). Our goal is to further understand the impact of the functional form of the bid on the market power of respective participants. In the case of the slope function bidding, each generator submits a slope function in the day-ahead and the real-time markets, parameterized by $\hat{b}_j^d \in \mathbb{R}_{\ge 0}, \ \hat{b}_j^r\in \mathbb{R}_{\ge 0}$, respectively: 
\begin{align}
     g_j^d = \hat{b}_j^d\lambda^d, \ \ g_j^r = \hat{b}_j^r\lambda^r 
\end{align}
Here $\lambda^d$ and $\lambda^r$ denote the prices in the day-ahead and real-time market, respectively. We first characterize the competitive equilibrium in a standard two-stage market. 

\begin{theorem}[Proposition  1~\cite{pyou_discovery}] \label{comp_eqbm_wout_mpm_slope} A competitive equilibrium in a two-stage market exists and is explicitly given by
\begin{subequations}
\begin{align}
    &{\hat{b}}_j^d+{\hat{b}}_j^r = \frac{1}{c_j}, \ {\hat{b}}_j^d \ge 0, \ {\hat{b}}_j^r \ge 0, \forall j \in \mathcal{G}\\
    & d_l^d+d_l^r = d_l, \forall l \in \mathcal{L}\\
    & \lambda^d = \lambda^r = \frac{d}{\sum_{j\in\mathcal{G}}c_j^{-1}}
\end{align}
\end{subequations}
\end{theorem}

The resulting competitive equilibrium is efficient, i.e., it aligns with the social planner problem~\eqref{planner_problem}. Similar to the competitive equilibrium for intercept function bidding in Theorem~\ref{comp_eqbm_wout_mpm_thrm}, the resulting equilibrium in Theorem~\ref{comp_eqbm_wout_mpm_slope} exists non-uniquely. We next consider the case of price-anticipating participants and characterize the resulting Nash equilibrium. 
\begin{theorem}[Proposition  4~\cite{pyou_discovery}] \label{strat_eqbm_wout_mpm_slope} Assume strategic generators are homogeneous $(c_j := c, \ \forall j \in  \mathcal{G})$. If there are at least three firms, i.e., $|\mathcal{G}| \ge 3$, a symmetric Nash equilibrium in a two-stage market exists with identical bids $(\hat{b}_j^v:= \hat{b}_j^v, \ \forall j \in  \mathcal{G},  v \in \{d,r\})$. Further, this equilibrium is unique and it is given by
\begin{align}
    &\!\!\!{\hat{b}}_j^d \!=\!\! \frac{|\mathcal{L}|(|\mathcal{G}|-1)+1}{|\mathcal{L}|(|\mathcal{G}|-1)}\frac{|\mathcal{G}|-2}{|\mathcal{G}|-1}\frac{1}{c}, \ \! {\hat{b}}_j^{r} \!\!=\!\! \frac{1}{|\mathcal{L}|+1}\frac{(|\mathcal{G}|-2)^2}{(|\mathcal{G}|-1)^2}\frac{1}{c} \label{strat_eqbm_wout_mpm_slope_eq.a}\\
    & \!\!\! d_l^{d} = \frac{|\mathcal{L}|(|\mathcal{G}|-1)+1}{|\mathcal{L}|(|\mathcal{L}|+1)(|\mathcal{G}|-1)}d, \ d_l^{r} \!\!= d_l - d_l^d \label{strat_eqbm_wout_mpm_slope_eq.b}\\
    & \lambda^{d} =  \frac{|\mathcal{L}|}{|\mathcal{L}|+1}\frac{|\mathcal{G}|-1}{|\mathcal{G}|-2}\frac{c}{|\mathcal{G}|}d, \ \lambda^{r} =  \frac{|\mathcal{G}|-1}{|\mathcal{G}|-2}\frac{c}{|\mathcal{G}|}d\label{strat_eqbm_wout_mpm_slope_eq.c}
\end{align}\label{strat_eqbm_wout_mpm_slope_eq}
\end{theorem}

\begin{figure}
  \centering
  \includegraphics[width=0.95\linewidth]{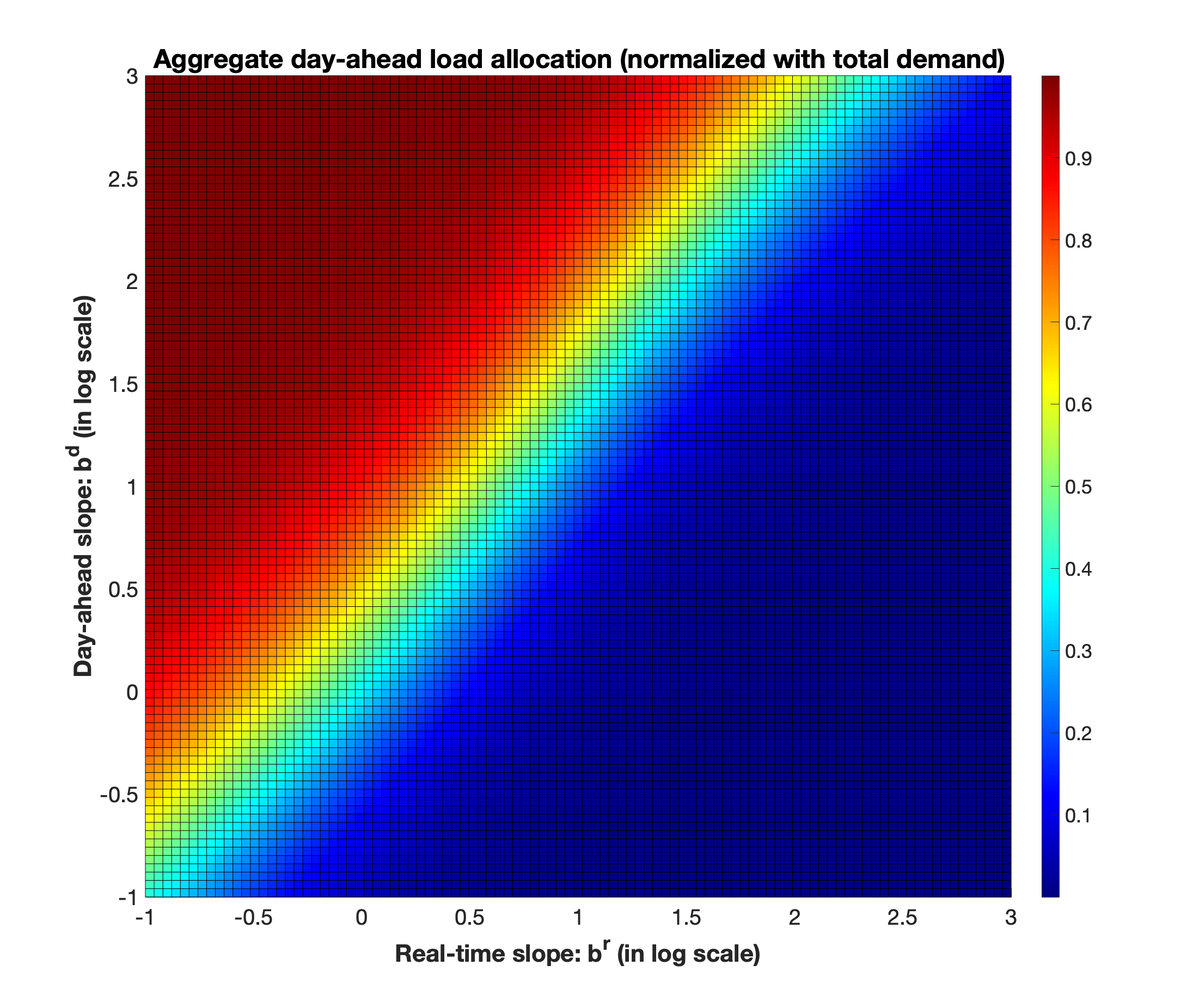}
  \caption{Normalized load allocation in the day-ahead stage in Intercept function bid-based standard market.}
  \label{fig:norm_da_allocation}
\end{figure}

Theorem \ref{strat_eqbm_wout_mpm_slope} shows the existence of a unique symmetric Nash equilibrium. At the resulting equilibrium, loads allocate more demand in the day-ahead market to exploit lower prices. However, the load allocation at the Nash equilibrium in the intercept function in Theorem~\ref{nash_eqbm_wout_mpm_thrm} is a function of market parameters $b^d$ and $b^r$. Figure~\ref{fig:norm_da_allocation} plots the aggregate load allocation in the day-ahead market as the slope of the intercept function bid changes in the day-ahead and real-time markets. We assume $4$ strategic homogeneous generators and $4$ strategic loads are participating in a standard two-stage market setting. The mix of individual inelastic demand bids is given by $d_l = [0.2, 25.6, 106.6, 199.6]^TMW$ from Pennsylvania,
New Jersey, and Maryland (PJM) data miner day-ahead demand bids~\cite{pjmdata} with total aggregate inelastic demand $d = 332 MW$. We assume a cost coefficient $c_j = 0.1 \$/MW^2,  \ \forall j \in \mathcal{G}$ corresponding to the cost coefficients from the IEEE 300-bus system~\cite{matpower} for homogeneous generators. The aggregate allocation in the day-ahead market (normalized with the total inelastic demand) can be increased by the operator with the help of appropriate slope parameters.     

\begin{figure*}
  \centering
  \includegraphics[width=0.9\linewidth]{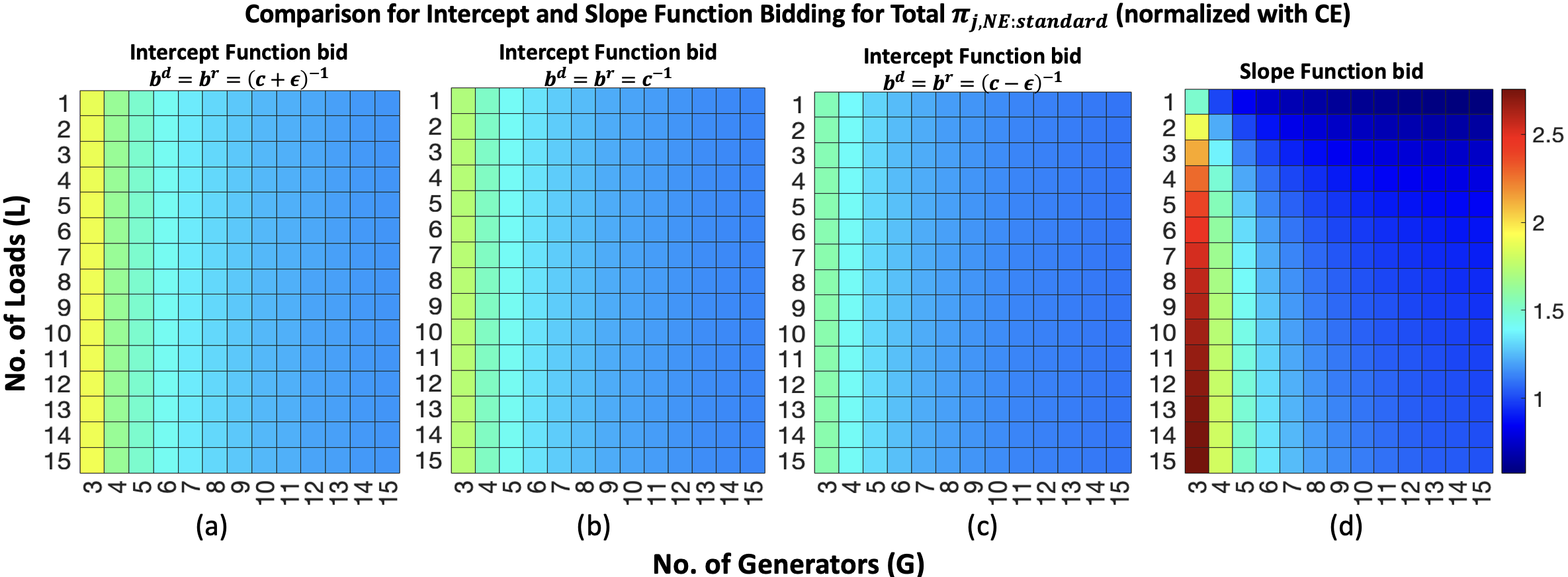}
  \caption{Aggregate generators' profit at Nash equilibrium (NE) normalized with competitive equilibrium (CE) in a standard market for Intercept function bid (a) with parameters $b^d = b^r = (c+\epsilon)^{-1}$, (b) with parameters $b^d = b^r = c^{-1}$, (c) with parameters $b^d = b^r = (c-\epsilon)^{-1}$, and (d) Slope function bid.}
  \label{fig:slope_intercept_compare}
\end{figure*}

Figure~\ref{fig:slope_intercept_compare} compares the (normalized) aggregate profit at Nash equilibrium in the standard market without any mitigation policy. We change the value of the slope parameter for the intercept function bid to understand the impact of model parameters, i.e., 
\[
    b^d = b^r = b, b \in \{(c+\epsilon)^{-1},c^{-1},(c-\epsilon)^{-1}\}, 
\]
where $\epsilon = 0.025 \$/MW^2$. The aggregate profit is normalized with the profit at competitive equilibrium levels. In the slope function bid-based market mechanism, there is a shift in the market power between loads and generators, e.g., loads win the competition for a relatively large number of generators in the market and vice versa. In particular, for a small number of loads and a large number of generators, loads exercise market power with lower payments at the expense of increased competition between generators. Similarly, a decrease in the number of generators and an increase in the number of loads favors generators in the market, as shown in panel (d) in Figure~\ref{fig:slope_intercept_compare}. However, generators always win the competition with higher profits at the Nash equilibrium in the intercept function bid based market mechanism, as shown in panel (b) in Figure~\ref{fig:slope_intercept_compare}. Moreover, such behavior, where generators always win the competition, exists regardless of slope parameter values in the intercept function bid, as shown in row 2 of Table~\ref{tab:cost_profit_table} and panels (a),(c) in Figure~\ref{fig:slope_intercept_compare}.

\section{Conclusions}\label{sec_6}

We model the competition between generators (bid intercept function) and loads (bid quantity) in a two-stage settlement market with stage-wise MPM policies. Per the CAISO's policy initiative,  non-competitive bids are substituted with the default bids based on estimated generator costs. We assume that a market with an MPM policy clears efficiently and generators participate truthfully in that stage, i.e., day-ahead or real-time. We start with a standard market without any mitigation policies. The competitive equilibrium in the standard market is efficient and non-unique, with an arbitrary allocation of demand across stages with equal prices. The Nash equilibrium results in a unique market outcome, and generators win the competition at the expense of loads with higher profits. 

Using the benchmark in the standard market setting, we first analyze a market with real-time MPM policy and show that the resulting competitive equilibrium is efficient. However, competition between loads and generators in the day-ahead market for individual interests drives all the demand to the real-time market. Finally, we analyze the day-ahead MPM policy market that results in efficient and unique competitive equilibrium. Moreover, our analysis of Nash equilibrium and its comparison with standard market outcome shows that the policy successfully mitigates the market power of generators and loads win the competition. 


\bibliographystyle{IEEEtran}
\bibliography{TSG}


\appendices
\section{Proof of Theorem~\ref{comp_eqbm_wout_mpm_thrm}}\label{appendix_comp_eqbm_wout_mpm_thrm}

Under price-taking behavior, the individual problem for loads~\eqref{load_price_taking_payment} is a linear program with the closed-form solution given by:
\begin{align}\label{comp_eqbm_load_solution_wout_mpm}
    \!\!\!\!\left\{\begin{array}{l}
d_l^{d} = \infty, d_l^{r} = -\infty, d_l^{d}+d_l^{r} = d_l, \!\mbox{ if } \lambda^{d} < \lambda^r \\
d_l^{d} = -\infty, d_l^{r} = \infty, d_l^{d}+d_l^{r} = d_l,  \mbox{ if } \lambda^{d} > \lambda^r \\
d_l^{d}+d_l^{r} = d_l, \quad  \mbox{ if } \lambda^{d} = \lambda^r
\end{array}\right.
\end{align}
where loads prefer the lower price in the market. The individual problem for generators~\eqref{gen_std_intercept_profit_comp} requires:
\begin{align}\label{comp_eqbm_gen_solution_wout_mpm}
    \!\!\!\!\!\!\left\{\begin{array}{l}
\!\beta_j^d \!= \!\infty, \beta_j^r \!=\! -\infty, \beta_j^d\!+\!\beta_j^r \!=\! \frac{b^d+b^r \!- c_j^{-1}}{\sum\limits_{j\in\mathcal{G}}\!\!c_j^{-1}}d, \!\mbox{ if } \lambda^{d} < \lambda^r \\
\!\beta_j^d \!=\! -\infty, \beta_j^r \!=\! \infty, \beta_j^d\!+\!\beta_j^r \!=\! \frac{b^d+b^r \!- c_j^{-1}}{\sum\limits_{j\in\mathcal{G}}\!\!c_j^{-1}}d, \!\mbox{ if } \lambda^{d} > \lambda^r \\
\!\beta_j^d+\beta_j^r = \frac{b^d+b^r - c_j^{-1}}{\sum\limits_{j\in\mathcal{G}}\!\!c_j^{-1}}d, \!\mbox{ if } \lambda^{d} = \lambda^r
\end{array}\right.
\end{align}
where generators prefer higher prices in the market and seek to maximize profit. 
At the competitive equilibrium the intercept function~\eqref{gen_intercept_bid_da},\eqref{gen_intercept_bid_rt} and individual optimal solution \eqref{comp_eqbm_load_solution_wout_mpm},\eqref{comp_eqbm_gen_solution_wout_mpm} holds simultaneously and this is only possible if the market price is equal in the two-stages, i.e., 
\[
    \lambda^{d} = \lambda^{r} = \frac{1}{\sum_kc_k^{-1}}d, d_l = d_l^d+d_l^r, \ \!\beta_j^d+\beta_j^r = \frac{b^d+b^r - c_j^{-1}}{\sum\limits_{j\in\mathcal{G}}\!\!c_j^{-1}}d
\]
From the intercept function bid, we have 
\[
    g_j^r+g_j^d = \frac{c_j^{-1}}{\sum_kc_k^{-1}}d
\]
Thus a set of competitive equilibria exists.

\section{Proof of Theorem~\ref{augmented_dispatch_thrm}} \label{strat_eqbm_wout_mpm_augmented_proof}

Given the parameter $(\beta_j^d, g_j^{d}, d-d^{d})$ from market-clearing in the day-ahead market, each generator $j$ maximizes their profit~\eqref{generator_strategic_profit_total} for the optimal decision $\beta_j^r$ with complete knowledge of the market clearing in the real-time stage as characterized below:
\begin{align}
    \!\!\!\!\!\!\!\!\sum_{j\in\mathcal{G}}g_j^{r} = d^{r} \!\!\!\! \implies \!\!\!\!\!\! \sum_{j\in\mathcal{G}}(b^r\lambda^{r} -\beta_j^r) = d^{r} 
    \!\!\!\!\implies \!\! \lambda^{r} = \frac{d^{r}+\beta^{r,\mathcal{G}}}{b^r|\mathcal{G}|} \label{gen_price_bid_function}
\end{align}
where $\beta^{r,\mathcal{G}} = \sum_{j\in\mathcal{G}}\beta_j^r$. Given the parameter $(\beta_j^d, g_j^{d}, d-d^{d})$, substituting~\eqref{gen_price_bid_function} in the individual problem~\eqref{generator_price_taking_profit} gives the concave strategic individual problem of generators, i.e., the real-time subgame problem: 
\label{augmented_dispatch_proof}
\begin{align}
    & \max_{\beta_j^r} \left(\frac{d^{r}+\beta^{r,\mathcal{G}}}{b^r|\mathcal{G}|}\right)\left(b^r\frac{d^{r}+\beta^{r,\mathcal{G}}}{b^r|\mathcal{G}|}-\beta_j^r\right)+\lambda^{d}g_j^{d} \nonumber \\
    & \quad \quad - \frac{c_j}{2}\left(g_j^{d}+b^r\left(\frac{d^{r}+\beta^{r,\mathcal{G}}}{b^r|\mathcal{G}|}\right)-\beta_j^r\right)^2 \label{gen_profit_startegic_thrm_augmented}
\end{align}
Hence, taking the derivative of~\eqref{gen_profit_startegic_thrm_augmented} with respect to bid $\beta_j^r$ we get:
\begin{align}
    &\frac{\partial \pi_j}{\partial b_j^r} = 0 \nonumber\\
    \implies & \frac{1}{b^r|\mathcal{G}|}\left(\frac{d^{r}+\beta^{r,\mathcal{G}}}{|\mathcal{G}|} - \beta_j^r\right) - \frac{|\mathcal{G}|-1}{|\mathcal{G}|}\left( \frac{d^{r}+\beta^{r,\mathcal{G}}}{b^r|\mathcal{G}|}\right)\nonumber \\
    & \quad +c_j\left(g_j^{d}+ \frac{d^{r}+\beta^{r,\mathcal{G}}}{|\mathcal{G}|}-\beta_j^r\right)\frac{|\mathcal{G}|-1}{|\mathcal{G}|} = 0 \nonumber\\
    \implies & \frac{1}{b^r|\mathcal{G}|}\left( g_j^{r}\right) - \frac{|\mathcal{G}|-1}{|\mathcal{G}|}\left( \lambda^{r}\right)+c_j\left(g_j^{d}+ g_j^{r}\right)\frac{|\mathcal{G}|-1}{|\mathcal{G}|} = 0 \nonumber\\
    \implies & \frac{1}{b(|\mathcal{G}|-1)}g_j^{r}-  \lambda^{r}+c_j\left(g_j^{d}+g_j^{r}\right) = 0 \label{augemented_obj_kkt_cond}
\end{align}
where we substitute~\eqref{gen_intercept_bid_rt} and~\eqref{gen_price_bid_function}. The equation~\eqref{augemented_obj_kkt_cond} is the required KKT condition of the convex dispatch problem~\eqref{augmented_dispatch}, with $\lambda^{r}$ as the dual variable of the constraint~\eqref{augmented_constraint}.

\section{Proof of Theorem~\ref{nash_eqbm_wout_mpm_thrm}} \label{app_nash_eqbm_wout_mpm}

Using the market-price in the real-time stage $\lambda^{r}$ as given by the KKT conditions~\eqref{augemented_obj_kkt_cond} we get, 
\begin{align}
    & g_j^{r} =  \frac{\lambda^{r}-c_jg_j^{d}}{C_j} \implies \sum_{j\in\mathcal{G}}g_j^{r} =\sum_{j\in\mathcal{G}}  \frac{\lambda^{r}-c_jg_j^{d}}{C_j} \nonumber \\
    \implies & \! d^{r} =\!\sum_{j\in\mathcal{G}}  \frac{\lambda^{r}-c_jg_j^{d}}{C_j} 
    \!\!\!\implies \!\!\! \lambda^{r} = \! \frac{d^{r} + \sum_{j\in\mathcal{G}} \frac{c_jg_j^{d}}{C_j}}{\sum_{j\in\mathcal{G}}C_j^{-1}} \label{gen_price_bid_function_startegic_thrm}
\end{align}
where $C_j =  \frac{1}{b^r(|\mathcal{G}|-1)}+c_j$ and we use~\eqref{rt_power_bal} in the second equality equation. Substituting~\eqref{gen_price_bid_function_startegic_thrm} in~\eqref{augemented_obj_kkt_cond} we get
\begin{align}
    g_j^{r} =  \frac{d^{r} + \sum_{k\in\mathcal{G}} \frac{c_kg_k^{d}}{C_k}}{C_j\sum_{k\in\mathcal{G}}      C_k^{-1} }-\frac{c_jg_j^{d}}{C_j} \label{gen_dispatch_bid_function_startegic_thrm}
\end{align}

From the market-clearing in the day-ahead stage~\eqref{da_power_bal}, we have the following relation 
\begin{subequations}
\begin{align}
     \implies &  \sum_{j \in \mathcal{G}} \left(b^d\lambda^d - \beta_j^d\right) = \sum_{l \in \mathcal{L}} d_l^d \\
    \implies & \! \lambda^{d} = \frac{d^{d}+\beta^{d,\mathcal{G}}}{b^d|\mathcal{G}|}, \ g_j^{d} = b^d\frac{d^{d}+\beta^{d,\mathcal{G}}}{b^d|\mathcal{G}|}-\beta_j^d \label{day_ahead_clearing_startegic_thrm}
\end{align}
\end{subequations}
where $\beta^{d,\mathcal{G}} = \sum_{j\in\mathcal{G}}\beta_j^d$. Substituting~\eqref{gen_price_bid_function_startegic_thrm},\eqref{gen_dispatch_bid_function_startegic_thrm},\eqref{day_ahead_clearing_startegic_thrm} in the individual profit~\eqref{generator_strategic_profit_total}, we get, 
\begin{small}
\begin{align}\label{generator_strategic_profit_total_wout_mpm}
    & \max_{\beta_j^d} \frac{d^d \!+\!\beta^{d,\mathcal{G}}}{b^d|\mathcal{G}|}\!\!\left(\!\frac{d^d \!+\! \beta^{d,\mathcal{G}}}{|\mathcal{G}|} - \beta_j^d\!\!\right) \!+\!\!\left(\!\frac{d^r \!+\!\!\! \sum\limits_{m\in\mathcal{G}} \!\!\frac{c_m}{C_m}\!\left(\!\frac{d^d \!+\! \beta^{d,\mathcal{G}}}{|\mathcal{G}|} \!- \beta_m^d\right)}{C_j\!\sum\limits_{k\in\mathcal{G}} C_k^{-1}}\!\right)^2 \nonumber \\
    &  -\frac{c_j}{C_j}  {\frac{d^r + \sum\limits_{m\in\mathcal{G}} \frac{c_m}{C_m}\left(\frac{d^d + \beta^{d,\mathcal{G}}}{|\mathcal{G}|} - \beta_m^d\right)}{\sum\limits_{k\in\mathcal{G}} C_k^{-1}}\left(\frac{d^d + \beta^{d,\mathcal{G}}}{|\mathcal{G}|} - b_j^d\right)}\nonumber \\
    & \!\!\!\!\!\!- \!\frac{c_j}{2}\!\!{\left(\!\!\!\left(\!\!1\!-\!\frac{c_j}{C_j}\!\!\right)\!\!\left(\!\frac{d^d \!+\! \beta^{d,\mathcal{G}}}{|\mathcal{G}|} \!- \beta_j^d\!\right)\!+\!\frac{d^r \!\!+\!\!\! \sum\limits_{m\in\mathcal{G}} \!\!\frac{c_m}{C_m}\!\!\left(\!\frac{d^d + \beta^{d,\mathcal{G}}}{|\mathcal{G}|} \!- \beta_m^d\!\right)}{C_j\!\sum\limits_{k\in\mathcal{G}} C_k^{-1}}\!\!\right)\!\!\!}^2 
\end{align}
\end{small}

Writing the first order condition and taking the derivative of~\eqref{generator_strategic_profit_total_wout_mpm} wrt $\beta_j^d$ we have
\begin{small}
\begin{align}\label{generator_strategic_profit_derivative_wout_mpm}
    & \implies  \frac{1}{b^d|\mathcal{G}|}\left(\!\frac{d^d \!+\! \beta^{d,\mathcal{G}}}{|\mathcal{G}|} - \beta_j^d\right) + \frac{d^d + \beta^{d,\mathcal{G}}}{b^d|\mathcal{G}|}\left(\frac{1}{|\mathcal{G}|} - 1\right)    \nonumber\\
    & + \frac{2}{C_j}\!\left(\frac{d^r + \sum\limits_{m\in\mathcal{G}} \frac{c_m}{C_m}\left(\frac{d^d + \beta^{d,\mathcal{G}}}{|\mathcal{G}|} - \beta_m^d\right)}{\sum_{k\in\mathcal{G}} C_k^{-1}}\right)\left(\frac{\sum\limits_{m\in\mathcal{G}} \frac{c_m}{C_m}\frac{1}{|\mathcal{G}|} - \frac{c_j}{C_j}}{\sum_{k\in\mathcal{G}} C_k^{-1}}\right)\nonumber \\
    & -\frac{c_j}{C_j}\!\!\left(\!\frac{\sum\limits_{m\in\mathcal{G}} \!\!\frac{c_m}{C_m}\frac{1}{|\mathcal{G}|} - \frac{c_j}{C_j}}{\sum_{k\in\mathcal{G}} C_k^{-1}}\!\right)\!\!\!\left(\frac{d^d \!+\! \beta^{d,\mathcal{G}}}{|\mathcal{G}|} - \beta_j^d\right) \nonumber\\
    & -\frac{c_j}{C_j}{\frac{d^r + \sum\limits_{m\in\mathcal{G}}\!\! \frac{c_m}{C_m}\left(\frac{d^d + \beta^{d,\mathcal{G}}}{|\mathcal{G}|} - \beta_m^d\right)}{\sum_{k\in\mathcal{G}} C_k^{-1}}\left(\frac{1}{|\mathcal{G}|} - 1\right)}\nonumber \\
    &\!\!\!\!\! - c_j\!\!\left(\!\!\!\!\left(1\!-\!\frac{c_j}{C_j}\right)\!\!\!\left(\frac{d^d + \beta^{d,\mathcal{G}}}{|\mathcal{G}|} - b_j^d\right)\!+\!\frac{d^r + \sum\limits_{m\in\mathcal{G}} \!\!\!\frac{c_m}{C_m}\!\!\left(\!\!\frac{d^d + \beta^{d,\mathcal{G}}}{|\mathcal{G}|} - b_m^d\right)}{C_j\sum\limits_{k\in\mathcal{G}} C_k^{-1}}\!\!\right) \nonumber \\ 
    & \quad \quad \quad\quad\quad \! \left(\!\!\!\left(1\!-\!\frac{c_j}{C_j}\right)\!\!\left(\frac{1}{|\mathcal{G}|} \!-\! 1\right)\!+\!\frac{1}{C_j} \frac{\sum\limits_{m\in\mathcal{G}}\!\! \frac{c_m}{C_m}\frac{1}{|\mathcal{G}|} - \frac{c_j}{C_j}}{\sum\limits_{k\in\mathcal{G}} C_k^{-1}}\!\right) \!=\! 0
\end{align}
\end{small}

Assuming generators are homogeneous, i.e. $c_j:= c, \ \forall j \in \mathcal{G}$ and we solve for symmetric equilibrium in the market, i.e., $\beta_j^d: = \beta^d, \ \forall j \in \mathcal{G}$, the equation \eqref{generator_strategic_profit_derivative_wout_mpm} can be rewritten as :
\begin{small}
\begin{align}\label{generator_strategic_profit_equilibrium_wout_mpm}
\implies &  \beta^d = b^dc\frac{d}{|\mathcal{G}|} +b^dc\frac{d^r}{|\mathcal{G}|}\left(1\!-\!\frac{c}{C}\right) \!-\! \frac{d^d}{|\mathcal{G}|}\frac{|\mathcal{G}|-2}{|\mathcal{G}|-1}
\end{align}
\end{small}
where $C :=  \frac{1}{b^r(|\mathcal{G}|-1)}+c$. 

Similarly, substituting~\eqref{gen_price_bid_function_startegic_thrm},\eqref{gen_dispatch_bid_function_startegic_thrm},\eqref{day_ahead_clearing_startegic_thrm} in the individual payment problem~\eqref{load_strategic_payment}, we get a convex optimization problem, 

\begin{align}\label{load_strategic_payment_wout_mpm}
    &\!\!\!\! \min_{d_l^d} \frac{d^d \!+\! \beta^{d,\mathcal{G}}}{b^d|\mathcal{G}|}d_l^d+\frac{d\!-\!d^d \!+\!\!\!\! \sum\limits_{m\in\mathcal{G}}\!\! \frac{c_m}{C_m}\!\!\left(\!\frac{d^d + \beta^{d,\mathcal{G}}}{|\mathcal{G}|} \!-\! \beta_m^d\!\right)}{\sum\limits_{k\in\mathcal{G}}C_k^{-1}}(d_l\!-\!d_l^d)
\end{align}

Taking the derivative of~\eqref{load_strategic_payment_wout_mpm} we have 
\begin{align}\label{load_strategic_payment_derivative_wout_mpm}
    \implies & \frac{d_l^d}{b^d|\mathcal{G}|} + \frac{d^d + \beta^{d,\mathcal{G}}}{b^d|\mathcal{G}|} +\frac{-1+ \sum\limits_{m\in\mathcal{G}}\!\! \frac{c_m}{C_m}\!\frac{1}{|\mathcal{G}|} }{\sum_{k\in\mathcal{G}}C_k^{-1}}(d_l-d_l^d) \nonumber \\
    & - \frac{d-d^d + \sum\limits_{m\in\mathcal{G}} \!\!\frac{c_m}{C_m}\!\!\left(\frac{d^d + \beta^{d,\mathcal{G}}}{|\mathcal{G}|} - \beta_m^d\right)}{\sum_{k\in\mathcal{G}}C_k^{-1}} = 0
\end{align}

Assume generators are homogeneous, i.e. $c_j:= c, \ \forall j \in \mathcal{G}$. We first sum over $l \in \mathcal{L}$ and solve for the case of symmetric bid participation of generators by rewriting the equation~\eqref{load_strategic_payment_derivative_wout_mpm} as,

\begin{small}
\begin{align} \label{load_strategic_payment_equilibrium_wout_mpm}
    \implies & d^d = - \frac{|\mathcal{G}|}{|\mathcal{L}|+1}\frac{|\mathcal{L}|\beta_j^d + b^dC\frac{-(|\mathcal{L}|+1)+ \frac{c}{C}}{|\mathcal{G}|}d}{1+\frac{b^d}{b^r(|\mathcal{G}|-1)}}
\end{align}
\end{small}

Solving the equations~\eqref{gen_intercept_bid_da},\eqref{gen_intercept_bid_rt},\eqref{gen_price_bid_function_startegic_thrm},\eqref{gen_dispatch_bid_function_startegic_thrm},\eqref{day_ahead_clearing_startegic_thrm},\eqref{generator_strategic_profit_equilibrium_wout_mpm}, and \eqref{load_strategic_payment_equilibrium_wout_mpm} simultaneously for the equilibrium, we get the unique Nash equilibrium as 
\begin{subequations}
    \begin{align}
          & \beta_j^d \!= \frac{b^dc}{|\mathcal{G}|}d + \frac{b^rc - \frac{|\mathcal{G}|-2}{|\mathcal{G}|-1}}{b^rc+\frac{|\mathcal{L}|+1}{|\mathcal{G}|-1}}\frac{|\mathcal{L}|+1}{|\mathcal{G}|(|\mathcal{G}|-1)}d^d, \ \forall j \in \mathcal{G} \\
          & \beta_j^r \!=\!  \frac{b^rc}{|\mathcal{G}|}d \!-\! \frac{|\mathcal{G}|-2}{|\mathcal{G}|(|\mathcal{G}|-1)}d^r, \ \forall j \in \mathcal{G} \\
           & \!d^{d} \!=\! \frac{b^d\left(b^rc+\frac{|\mathcal{L}|+1}{|\mathcal{G}|-1}\right)}{b^d\!\left(b^rc\!+\!\frac{|\mathcal{L}|+1}{|\mathcal{G}|-1}\right) \!+\! b^r\!\left(b^rc\!+\!\frac{1}{|\mathcal{G}|-1}\!\right)\!(|\mathcal{G}|+|\mathcal{L}|\!-\!1)\!}d \\
          & \!\!d_l^d \!\!= \!\frac{b^dd_l}{b^d\!+\!b^r\!(\!|\mathcal{G}|\!\!-\!1\!)} \!+\! \frac{\frac{b^d}{1+b^rc(|\mathcal{G}|\!-\!1)}}{b^d\!+\!b^r\!(\!|\mathcal{G}|\!\!-\!1\!)}d^r \!\!\!-\! \frac{b^r}{b^d\!+\!b^r\!(\!|\mathcal{G}|\!\!-\!1\!)}d^d\\
          & \lambda^d \!= \!\frac{b^rc(|\mathcal{G}|\!-\!1) \!+\!2}{b^rc(|\mathcal{G}|\!-\!1) \!+\! 1}\frac{c}{|\mathcal{G}|}d+\frac{\left(\frac{b^r}{b^d}\!-\!1\right)c \!+\! \frac{1}{b^d(|\mathcal{G}|-1)}}{b^rc(|\mathcal{G}|\!-\!1) \!+\! 1}\frac{d^d}{|\mathcal{G}|}, \\
          & \!\!\lambda^r \!\!=\! \lambda^d \!+\!\frac{\frac{1}{|\mathcal{G}|(|\mathcal{G}|\!-\!1)}\left(\frac{|\mathcal{G}|\!-\!2}{|\mathcal{G}|\!-\!1} \!-\! b^rc\right)d}{b^d\!\left(\!b^rc\!+\!\frac{|\mathcal{L}|+1}{|\mathcal{G}|-1}\!\right) \!\!+\! b^r\!\!\left(\!b^rc\!+\!\frac{1}{|\mathcal{G}|-1}\!\right)\!(|\mathcal{G}|\!+\!|\mathcal{L}|\!\!-\!\!1\!)\!} 
    \end{align}
    \end{subequations}

Thus the symmetric Nash equilibrium exists uniquely. 

\section{Proof of Theorem~\ref{comp_eqbm_rt_mpm_thrm}} \label{app_comp_eqbm_rt_mpm}

Under price-taking behavior, the individual problem for loads~\eqref{load_payment_intermediate_rt_mpm} is a linear program with the closed-form solution given by:
\begin{align}\label{comp_eqbm_load_solution_rt_mpm}
    \!\!\!\!\left\{\begin{array}{l}
d_l^{d} = \infty, d_l^{r} = -\infty, d_l^{d}+d_l^{r} = d_l, \!\mbox{ if } \lambda^{d} < \frac{d}{\sum_kc_k^{-1}} \\
d_l^{d} = -\infty, d_l^{r} = \infty, d_l^{d}+d_l^{r} = d_l,  \mbox{ if } \lambda^{d} > \frac{d}{\sum_kc_k^{-1}}  \\
d_l^{d}+d_l^{r} = d_l, \quad  \mbox{ if } \lambda^{d} = \frac{d}{\sum\limits_{k \in \mathcal{G}}c_k^{-1}}
\end{array}\right.
\end{align}
where loads prefer the lower price in the market. The individual problem for generators~\eqref{generator_profit_rt_mpm} requires:
\begin{align}\label{comp_eqbm_gen_solution_rt_mpm}
    \left\{\begin{array}{l}
\beta_j^d = \infty, \!\mbox{ if } \lambda^{d} < \frac{d}{\sum\limits_{k \in \mathcal{G}}c_k^{-1}} \\
\beta_j^d = -\infty,  \mbox{ if } \lambda^{d} > \frac{d}{\sum\limits_{k \in \mathcal{G}}c_k^{-1}}  \\
\beta_j^d \in \mathbb{R}, \quad  \mbox{ if } \lambda^{d} = \frac{d}{\sum\limits_{k \in \mathcal{G}}c_k^{-1}}
\end{array}\right.
\end{align}
where generators prefer higher prices in the market and seek to maximize profit. At the competitive equilibrium the day-ahead supply function~\eqref{gen_intercept_bid_da}, real-time true dispatch condition~\eqref{rt_true_dispatch}, real-time clearing prices~\eqref{rt_true_prc}, and the individual optimal solution \eqref{comp_eqbm_load_solution_rt_mpm},\eqref{comp_eqbm_gen_solution_rt_mpm} holds simultaneously and this is only possible if the market price is equal in the two-stages, i.e., 
\[
    \lambda^{d} = \lambda^{r} = \frac{1}{\sum\limits_{k \in \mathcal{G}}c_k^{-1}}d, \textrm{s.t } d_l = d_l^d+d_l^r, \ \beta_j^d \in \mathbb{R}
\]
and 
\[
    g_j^r+g_j^d = \frac{c_j^{-1}}{\sum_kc_k^{-1}}d
\]
Thus a set of competitive equilibria exists.

\section{Proof of Theorem~\ref{strat_eqbm_rt_mpm_thrm}}\label{app_strat_eqbm_rt_mpm}

Substituting the real-time true dispatch condition~\eqref{rt_true_dispatch}, real-time clearing prices~\eqref{rt_true_prc}, day-ahead dispatch and day-ahead prices~\eqref{day_ahead_clearing_startegic_thrm} in the individual problem of generator~\eqref{generator_strategic_profit_rt_mpm}, we get

\begin{small}
\begin{align}\label{gen_strat_profit_rt_mpm_proof}
    &\!\!\!\!\!\! \max_{\beta_j^d}\!\!\! \ \left(\!\!\frac{d^{d}\!+\!\beta^{d,\mathcal{G}}}{b^d|\mathcal{G}|}-\frac{d}{\sum\limits_{k \in  \mathcal{G}}c_k^{-1}}\!\!\!\right)\!\!\!\left(\!\frac{d^{d}\!+\!\beta^{d,\mathcal{G}}}{|\mathcal{G}|}\!-\beta_j^d\!\!\right)\! +\!\frac{c_j^{-1}}{2}\!{\left(\!\frac{d}{\sum\limits_{j \in  \mathcal{G}}c_j^{-1}}\!\!\right)\!\!}^2    
\end{align}
\end{small}
where $\beta^{d,\mathcal{G}} = \sum_{j\in\mathcal{G}}\beta_j^d$. Taking the derivative of~\eqref{gen_strat_profit_rt_mpm_proof} wrt $\beta_j^d$ and writing the first-order condition, we have 
\begin{small}
\begin{align}\label{gen_strat_deriv_rt_mpm_proof}
    & \!\!\!\!\!\! \frac{1}{b^d|\mathcal{G}|}\!\left(\!\frac{d^{d}+\beta^{d,\mathcal{G}}}{|\mathcal{G}|}-\beta_j^d\!\right) \!\!+\!\! \left(\!\!\frac{d^{d}+\beta^{d,\mathcal{G}}}{b^d|\mathcal{G}|}-\frac{d}{\sum\limits_{k \in  \mathcal{G}}c_k^{-1}}\!\!\right)\!\!\!\left(\!\frac{1}{|\mathcal{G}|}\!-\!1\!\right) \!=\! 0
\end{align}
\end{small}
Summing the equation~\eqref{gen_strat_deriv_rt_mpm_proof} over the set of generators, i.e., $j \in \mathcal{G}$ we get
\begin{small}
\begin{subequations}
\begin{align}\label{gen_strat_deriv_contd_rt_mpm_proof}
    \implies & \frac{1}{b^d|\mathcal{G}|}d^{d} - \left(\frac{d^{d}+\beta^{d,\mathcal{G}}}{b^d|\mathcal{G}|}-\frac{d}{\sum\limits_{k \in  \mathcal{G}}c_k^{-1}}\right)\left(|\mathcal{G}|-1\right) \!=\! 0 \\
    \implies & \beta^{d,\mathcal{G}} = \frac{b^d|\mathcal{G}|}{\sum\limits_{k \in  \mathcal{G}}c_k^{-1}}d - \frac{(|\mathcal{G}|-2)}{(|\mathcal{G}|-1)}d^{d} \label{gen_strat_deriv_contd_rt_mpm_proof.a}\\
    \implies & \beta_j =  b^d\frac{d}{\sum\limits_{k \in  \mathcal{G}}c_k^{-1}} - \frac{|\mathcal{G}|-2}{|\mathcal{G}|} \frac{1}{(|\mathcal{G}|-1)}d^{d}\label{gen_strat_deriv_contd_rt_mpm_proof.b}
\end{align}
\end{subequations}
\end{small}

Similarly, substituting the real-time clearing prices~\eqref{rt_true_prc} and day-ahead prices~\eqref{day_ahead_clearing_startegic_thrm} in the individual problem of generator~\eqref{load_strategic_payment_rt_mpm}, we get

\begin{small}
\begin{align}\label{load_strat_pay_rt_mpm_proof}
    & \min_{d_l^{d}} \ \left(\frac{d^{d}+\beta^{d,\mathcal{G}}}{b^d|\mathcal{G}|}\right)d_l^{d} +\left(\frac{d}{\sum\limits_{k \in  \mathcal{G}}c_k^{-1}}\right)(d_l-d_l^{d})    
\end{align}
\end{small}

Taking the derivative of ~\eqref{load_strat_pay_rt_mpm_proof} and writing the first order condition of the convex optimization problem, we get
\begin{align}\label{load_strat_pay_deriv_rt_mpm_proof}
    \frac{d_l^{d}+ d^{d}+\beta^{d,\mathcal{G}}}{b^d|\mathcal{G}|} - \frac{d}{\sum\limits_{k \in  \mathcal{G}}c_k^{-1}} = 0
\end{align}
Summing the equation~\eqref{load_strat_pay_deriv_rt_mpm_proof} over $l\in\mathcal{L}$, we get
\begin{align}\label{load_strat_pay_deriv_rt_mpm_proof.a}
    \implies d^{d} = \frac{|\mathcal{L}|}{|\mathcal{L}|+1}\frac{b^d|\mathcal{G}|}{\sum\limits_{k \in  \mathcal{G}}c_k^{-1}}d - \frac{|\mathcal{L}|}{|\mathcal{L}|+1}\beta^{d,\mathcal{G}}
\end{align}

At the equilibrium the equations~\eqref{load_two_stage},\eqref{rt_true_prc},\eqref{day_ahead_clearing_startegic_thrm},\eqref{gen_strat_deriv_contd_rt_mpm_proof.b}, and \eqref{load_strat_pay_deriv_rt_mpm_proof.a} must hold simultaneously. Solving for equilibrium we get, 
\begin{small}
\begin{subequations}
\begin{align}     
& d_l^{d} = 0, \ d_l^{r} = d_l, \ \forall l\in\mathcal{L} \\ 
& \lambda^{d} = \lambda^{r} = \frac{1}{\sum\limits_{j\in\mathcal{G}}c_j^{-1}}d   \\
 & g_j^{d} = 0,  \ g_j^{r} = \frac{c_j^{-1}}{\sum\limits_{j\in\mathcal{G}}c_j^{-1}}d, \ \beta_j^{d} = \frac{b^d}{\sum\limits_{j \in  \mathcal{G}}c_j^{-1}}d, \ \!\forall j\in\mathcal{G} 
\end{align}\label{strat_eqbm_rt_mpm}
\end{subequations}
\end{small}

This completes the proof.

\section{Proof of Theorem~\ref{comp_eqbm_da_mpm}}\label{app_comp_eqbm_da_mpm}

Under price-taking behaviour, the individual problem for loads~\eqref{load_price_taking_payment} is a linear program with the closed-form solution given by:
\begin{align}\label{comp_eqbm_load_solution}
    \left\{\begin{array}{l}
d_l^{d} = \infty, d_l^{r} = -\infty, d_l^{d}+d_l^{r} = d_l, \mbox{ if } \lambda^{d} < \lambda^{r}  \\
d_l^{d} = -\infty, d_l^{r} = \infty, d_l^{d}+d_l^{r} = d_l,  \mbox{ if } \lambda^{d} > \lambda^{r}  \\
d_l^{d}+d_l^{r} = d_l, \quad  \mbox{ if } \lambda^{d} = \lambda^{r}
\end{array}\right.
\end{align}
where loads prefer the lower price in the market. Further solving the individual bidding problem for generators in real-time market~\eqref{generator_price_taking_profit_bids_da_mpm} by taking the derivative of the concave profit function wrt $\beta_j^r$, we get

\begin{align}
    &-\lambda^{r}+c_j\left(\frac{c_j^{-1}d^d}{\sum\limits_{k\in\mathcal{G}}c_k^{-1}}+b^r\lambda^{r}-\beta_j^r\right) = 0 \label{comp_eqbm_generator_solution}
\end{align}
Substituting~\eqref{da_true_dispatch},\eqref{da_true_prc} and \eqref{gen_intercept_bid_rt}
in~\eqref{comp_eqbm_generator_solution}, we get 
\begin{align}
    \implies & -\lambda^{r}+c_j(g_j^{d}+g_j^{r}) = 0
    \implies \sum_{j\in\mathcal{G}}\frac{1}{c_j}\lambda^{r} = \sum_{j\in\mathcal{G}}g_j = d \nonumber \\
    \implies & \lambda^{r} = \frac{d}{\sum\limits_{j\in\mathcal{G}}c_j^{-1}} \label{comp_eqbm_gen_real_time_price}
\end{align}
At the competitive equilibrium the conditions~\eqref{da_true_dispatch},\eqref{da_true_prc},\eqref{comp_eqbm_load_solution},\eqref{comp_eqbm_generator_solution}, and \eqref{comp_eqbm_gen_real_time_price} must hold simultaneously and this is only possible if the market price are equal in the two-stages, i.e., 
\[
    \lambda^{r} = \lambda^{d} = \frac{d}{\sum_{j\in\mathcal{G}}c_j^{-1}}
\]
and 
\[
    d^{d} = d, \ d^{r} = 0; \  d_l^{d}+d_l^{r} = d_l,  \forall l\in\mathcal{L},
\]
\[
    g_j^{d} = \frac{1}{c_j}\frac{d}{\sum\limits_{k\in\mathcal{G}}c_k^{-1}}, \ g_j^{r} = 0, \ \forall j\in\mathcal{G}
\]
Thus the competitive equilibrium exists.

\section{Proof of Theorem~\ref{strat_eqbm_da_mpm_thrm}}\label{app_strat_eqbm_da_mpm}

Using the Theorem $3.4$ in reference~\cite{bansal_e_energy} and the market-price in the real-time stage $\lambda^{r}$ as given by the KKT 
conditions~\eqref{augemented_obj_kkt_cond}, we get 
\begin{align}
    & \lambda^{r} = \! \frac{d^{r} + \sum_{j\in\mathcal{G}} \frac{c_jg_j^{d}}{C_j}}{\sum_{j\in\mathcal{G}}C_j^{-1}} \label{gen_price_bid_function_startegic_thrm_2}\\
    & g_j^{r} =  \frac{d^{r} + \sum_{j\in\mathcal{G}} \frac{c_jg_j^{d}}{C_j}}{\sum_{j\in\mathcal{G}}      C_j^{-1} C_j}-\frac{c_jg_j^{d}}{C_j} \label{gen_dispatch_bid_function_startegic_thrm_2}
\end{align}
where $C_j =  \frac{1}{b^r(|\mathcal{G}|-1)}+c_j$.
Substituting~\eqref{da_true_dispatch},\eqref{da_true_prc} in the expression~\eqref{gen_price_bid_function_startegic_thrm_2} and~\eqref{gen_dispatch_bid_function_startegic_thrm_2} we get
\begin{align}
    & \lambda^{r} = \frac{d^{r}}{\sum\limits_{k\in\mathcal{G}} C_k^{-1}}+ \frac{d^{d}}{\sum\limits_{k\in\mathcal{G}} c_k^{-1}}, \ g_j^{r} =  \frac{1}{C_j}\frac{d^{r}}{\sum\limits_{k\in\mathcal{G}} C_k^{-1}} \label{price_dispatch_real_time_startegic_thrm}
\end{align}
Substituting~\eqref{da_true_prc} and~\eqref{price_dispatch_real_time_startegic_thrm} in the individual problem of load $l$~\eqref{load_strategic_payment} we get
\begin{align}
    & \min_{d_l^{d}} \ \!\frac{d^{d}}{\sum_{j \in G}c_j^{-1}}d_l^{d} \!+\! \left(\!\frac{d-d^{d}}{\sum\limits_{k\in\mathcal{G}}C_k^{-1}}\!+\! \frac{d^{d}}{\sum\limits_{k\in\mathcal{G}}c_k^{-1}}\!\right)\!\!(d_l-d_l^{d}) \label{load_payment_startegic_thrm}
\end{align}
Therefore taking the derivative of the convex individual problem~\eqref{load_payment_startegic_thrm} wrt $d_l^{d}$ we get,
\begin{align}
    & - \frac{d-d^{d}}{\sum\limits_{k\in\mathcal{G}}      C_k^{-1}} + \frac{d_l}{\sum\limits_{k\in\mathcal{G}} c_k^{-1}}-\frac{d_l}{\sum\limits_{k\in\mathcal{G}}C_k^{-1}} + \frac{d_l^{d}}{\sum\limits_{k\in\mathcal{G}}      C_k^{-1}} = 0 \nonumber \\
    \implies & \!\!\!\! \sum_{l\in\mathcal{L}}\!\!\left(\!- \frac{d-d^{d}}{\sum\limits_{k\in\mathcal{G}}C_k^{-1}} + \frac{d_l}{\sum\limits_{k\in\mathcal{G}}c_k^{-1}}-\frac{d_l}{\sum\limits_{k\in\mathcal{G}}C_k^{-1}} + \frac{d_l^{d}}{\sum\limits_{k\in\mathcal{G}}C_k^{-1}}\!\!\right) \!=\! 0 \nonumber \\
    \implies & d^{d} = \left(1- \frac{1}{|\mathcal{L}|+1}\frac{\sum\limits_{k\in\mathcal{G}}C_k^{-1}}{\sum\limits_{k\in\mathcal{G}}c_k^{-1}}\right)d \label{day_ahead_dispatch_strategic_theorem}
\end{align}
Therefore we get unique Nash equilibrium ~\eqref{strat_eqbm_traditional_equations}

\end{document}